\PassOptionsToPackage{cmex10}{amsmath}
\documentclass[12pt]{amsart}
\usepackage{geometry}
\usepackage{amsmath,amssymb,amsthm,mathtools}
\usepackage{enumitem}
\usepackage{xcolor}
\usepackage[colorlinks,citecolor=red,linkcolor=blue,urlcolor=blue,hyperindex]{hyperref}
\usepackage{bookmark}
\DeclareFontShape{OT1}{cmr}{bx}{n}{<-8>cmbx7<8-10>cmbx8<10->cmbx10}{}
\DeclareFontShape{OMX}{cmex}{m}{n}{<->cmex10}{}
\allowdisplaybreaks[2]
\newtheorem{theorem}{Theorem}[section]
\newtheorem{lemma}[theorem]{Lemma}
\newtheorem{proposition}[theorem]{Proposition}
\newtheorem{corollary}[theorem]{Corollary}
\theoremstyle{definition}

\theoremstyle{remark}
\newtheorem{remark}{Remark}[section]
\numberwithin{equation}{section}
\newcommand{\D}{\mathbb D}

\newcommand{\norm}[1]{\left\|#1\right\|}
\newcommand{\dd}{\,\mathrm{d}}
\newcommand{\z}{\mathbf z}
\newcommand{\w}{\mathbf w}
\newcommand{\aA}{\mathbf a}
\newcommand{\m}{\mathbf m}
\newcommand{\om}{\boldsymbol\omega}
\newcommand{\la}{\boldsymbol\lambda}
\newcommand{\Hb}{h_{\bar b}^{\aA}}
\newcommand{\T}{T_{\aA,\m}}
\newcommand{\Bp}{\mathcal B_\Pi}
\newcommand{\pair}[2]{\langle #1,#2\rangle_{\aA}}
\newcommand{\bd}[2]{(1-|#1|^2)^{#2}}
\newcommand{\ow}[1]{\omega(1-|#1|)}
\begin{document}
\title[Small Hankel operators on weighted Besov spaces]{Small Hankel Operators
on Weighted Besov Spaces over the Polydisc}
\author{Chen Li}
\address{School of Mathematical Sciences, Fudan University, Shanghai, China}
\email{22110180024@m.fudan.edu.cn}
\author{Kai Wang}
\address{School of Mathematical Sciences, Fudan University, Shanghai, China}
\email{kwang@fudan.edu.cn}
\date{}
\begin{abstract}
We characterize bounded small Hankel operators from $B_p(\om)$ to
$\overline{B_q(\om)}$ on the polydisc for $0<p,q<\infty$, under
explicit conditions relating the kernel parameter to the weight.
When $q<p$, boundedness is equivalent to a weighted Besov condition
on the symbol; when $p\le q$, it is equivalent to a weighted
supremum condition on a shifted radial derivative.
In particular, the diagonal case is characterized by the product
Bloch space. The proof combines weighted Bergman duality and atomic
factorization with local kernel estimates and randomized atomic tests.
\end{abstract}
\maketitle

\section{Introduction}\label{sec:main}

We characterize bounded small Hankel operators between weighted
holomorphic Besov spaces on the polydisc for $0<p,q<\infty$.
The source and target spaces carry the same product weight, while
the operator is defined by a standard weighted Bergman kernel.
Its parameters are chosen independently of the space weight, subject
to the conditions stated below.

Symbol characterizations of Hankel operators are closely related to
duality and factorization in Bergman spaces. For Hilbert Bergman spaces
on the disc with different power weights, Janson \cite{Janson} obtained
Besov-space criteria for boundedness, compactness and Schatten-class
membership of both big and small Hankel operators. On bounded symmetric
domains, Zhu \cite{Zhu} characterized these properties for little Hankel
operators by integral transforms of their symbols, with further results
for the ball and the polydisc. Bonami and Luo \cite{BonamiLuo}
studied small Hankel operators between unweighted Bergman spaces
on the unit ball for positive source and target exponents, including
exponents below one. Pau and Zhao \cite{PauZhao} studied bounded
Hankel forms and weak factorization on power-weighted Bergman spaces of
the unit ball. Their small Hankel results allow different power weights
and compatible kernel parameters; in particular, for $1<q<p<\infty$,
the same-weight case is characterized by a Bergman integrability condition
with exponent $pq/(p-q)$.

For radial weights on the disc satisfying a doubling condition,
Korhonen and R\"atty\"a
\cite[Corollary~6]{KorhonenRattya} used factorization to obtain the corresponding
small Hankel characterization in the range $1<q<p<\infty$.
Duan, R\"atty\"a, Wang and Wu \cite{Wang} treated the full range
$0<p,q<\infty$ with two radial weights satisfying two-sided doubling
conditions. In their setting the projection is induced by the target
weight. Pel\'aez and R\"atty\"a \cite[Theorem~1]{PelaezRattya} also obtained a
criterion for measurable symbols acting on a single Bergman space
with a two-sided doubling radial weight and exponent $1<p<\infty$,
using modified Bergman kernels.

A related theory concerns big Hankel operators, whose values lie in the
complement of the holomorphic subspace. For general symbols on the
disc, Luecking \cite{LueckingHankel} characterized boundedness from
$A^p$ to $L^p$, $1<p<\infty$, and studied compactness and
Schatten-class membership in the Hilbert-space case. Galanopoulos
and Pau \cite{GalanopoulosPau} obtained derivative criteria for
boundedness and compactness with conjugate-analytic symbols on
Hilbert Bergman spaces with rapidly decreasing radial weights.
On the unit ball, Pau, Zhao and
Zhu \cite{PauZhaoZhu} characterized simultaneous boundedness and
compactness of the operators with symbols $f$ and $\overline f$ by
weighted BMO and VMO conditions, respectively, for $1<p\le q<\infty$.
Pel\'aez, Per\"al\"a and R\"atty\"a \cite[Theorem~1]{PelaezPeralaRattya}
proved a weighted mean-oscillation characterization of simultaneous
boundedness on the disc in the same exponent range, with a two-sided doubling
source weight and a radial B\'ekoll\'e--Bonami target weight.
For the small Hankel operators considered here, the symbol
conditions take the form of analytic integrability or derivative growth.

The use of kernel atoms in representation and factorization theorems
goes back to Coifman and Rochberg \cite{CoifmanRochberg}.
On the polydisc, Constantin's atomic decompositions \cite{Constantin}
and weak product results \cite{ConstantinWeak} provide a way to relate
Hankel forms to symbol functionals, including in vector-valued weighted
Bergman spaces. On the ball, the connection between weak
factorization and Hankel forms extends to weighted Bergman--Orlicz
spaces \cite{TchoundjaZhao,BekolleBonamiTchoundja}; the latter work
also uses atomic decompositions with concave growth functions and
factorizations involving a Bloch factor.

For the Besov spaces considered here, Harutyunyan and
Marinescu \cite{HM} studied generalized little Hankel and Berezin-type
operators with bounded measurable symbols. Their results concern equal
source and target exponents and include exponents below one, subject
to restrictions connecting the power-kernel parameters with the weight
indices. Here we allow the source and target Besov exponents to differ
and characterize boundedness in terms of the symbol. For the standard
small Hankel operator on the unweighted polydisc, Sehba
\cite[Theorem~2.1]{Sehba} identifies a logarithmic Bloch condition
at the endpoint $A^1$. That matched-kernel endpoint is excluded
by our strict kernel assumptions; see the remark after
Proposition~\ref{prop:outputdual}.

Our characterization has two forms: an integral condition when $q<p$
and a weighted derivative growth condition when $p\le q$.
The proof follows the atomic approach to Hankel forms when $q\ge1$.
When $q<1$, local output estimates and randomized tests replace the
duality argument. We begin with the definitions.

Let $\D^n=\{\z=(z_1,\ldots,z_n):|z_j|<1\}$ and
$dV(\z)=\prod_{j=1}^n dA(z_j)$, where $dA=\pi^{-1}dx\,dy$.
We write $H(\D^n)$ for the holomorphic functions on $\D^n$ and use bold
letters for vectors and multi-indices. The notation $X\lesssim Y$ means
$X\le CY$ with a constant independent of the functions and sampling
points; $X\asymp Y$ means that the estimate holds in both directions.

Motivated by the weighted Besov spaces of Harutyunyan and Lusky
\cite{HL}, we consider, in each coordinate, a weight of the form
\[
 \omega_j(t)=\exp\left(\eta_j(t)+\int_t^1\varepsilon_j(u)\frac{du}{u}\right),
 \qquad 0<t\le1,
\]
where $\eta_j$ is bounded and
$-\gamma_j^+\le\varepsilon_j\le-\gamma_j^-$, with
$1<\gamma_j^-\le\gamma_j^+$. These assumptions imply
\begin{equation}\label{eq:scale}
 C_j^{-1}(s/t)^{\gamma_j^+}
 \le\frac{\omega_j(s)}{\omega_j(t)}
 \le C_j(s/t)^{\gamma_j^-},\qquad 0<s\le t\le1.
\end{equation}
Only the two-sided estimate \eqref{eq:scale} will be used below, so the
argument also applies to positive measurable weights satisfying it.
Set $\om=(\omega_1,\ldots,\omega_n)$ and
$R_j=z_j\frac{\partial}{\partial z_j},\, \mathfrak D=\prod_{j=1}^n(I+R_j).$
For $0<t<\infty$, the weighted Besov space is defined by the condition
\begin{equation}\label{eq:besov}
 \norm f_{B_t(\om)}^t
 =\int_{\D^n}|\mathfrak D f(\z)|^t
 \prod_j\omega_j(1-|z_j|)(1-|z_j|^2)^{t-2}\dd V(\z)<\infty.
\end{equation}
The anti-holomorphic space is defined by conjugation. For $t<1$, we
use the corresponding quasi-norm.

To specify the kernel, let $\aA=(a_1,\ldots,a_n)$, $a_j>-1$, and set
\[
 d\nu_{\aA}(\w)=\prod_j(a_j+1)(1-|w_j|^2)^{a_j}\dd A(w_j),
 \qquad
 K_{\aA}(\w,\z)=\prod_j(1-w_j\overline{z_j})^{-a_j-2}.
\]
Write $A_{\aA}^t=H(\D^n)\cap L^t(\nu_{\aA})$ for the
standard weighted Bergman space.
For $b\in H(\D^n)\cap L^1(\nu_{\aA})$, the small Hankel operator is
\begin{equation}\label{eq:hankel}
 (\Hb f)(\z)=\int_{\D^n}f(\w)\overline{b(\w)}
 K_{\aA}(\w,\z)\dd\nu_{\aA}(\w).
\end{equation}
This integral defines an anti-holomorphic function for every input
in $\mathcal H_{\mathrm{cl}}$, the class of functions holomorphic
in a neighbourhood of the closed polydisc. Boundedness means that
this operator extends continuously to the stated space. For general
inputs, we use the shifted integral representation whenever its
absolute convergence has been established.

We write $q'=q/(q-1)$ for $q>1$ and set $q'=\infty$ at $q=1$,
with $1/\infty=0$.

For products of coordinate weights and boundary factors, we use the notation
$\omega(\mathbf t)=\prod_{j=1}^n\omega_j(t_j),\, \bd{\z}{\mathbf c}=\prod_{j=1}^n(1-|z_j|^2)^{c_j}.$
If the exponent is a scalar $c$, each $c_j$ in this convention is
$c$; sums such as $\aA+c$ are also interpreted coordinatewise.
Thus $\om=(\omega_1,\ldots,\omega_n)$ denotes the tuple of weights,
whereas $\omega(\mathbf t)$ denotes the product of its components. The density
associated with the Bergman-space formulation of the problem is
\begin{equation}\label{eq:weights}
 v_j(x)=\frac{\omega_j(1-x)}{(1-x^2)^2},\qquad
 v(\z)=\prod_jv_j(|z_j|)
       =\frac{\ow{\z}}{\bd{\z}{2}}.
\end{equation}
The symbol conditions will be expressed through a differential operator
adapted to the kernel. For positive integers $m_j$, let
\begin{equation}\label{eq:T}
 T_{\aA,\m}=\prod_{j=1}^n\prod_{\ell=0}^{m_j-1}
       \frac{R_j+a_j+2+\ell}{a_j+2+\ell}.
\end{equation}
The lower-order terms ensure that $T_{\aA,\m}1=1$ and allow the
operator to detect constant symbols. Lemma~\ref{lem:shift} explains
how this derivative arises from a shift of the weight in the pairing.

\begin{theorem}\label{thm:main}
Let $0<p,q<\infty$ and $a_j>-1$, $j=1,\ldots,n$.
Suppose that the positive weights $\omega_j$ satisfy \eqref{eq:scale}
with $1<\gamma_j^-\le\gamma_j^+$, and assume
\begin{equation}\label{eq:kernelcondition}
 \begin{cases}
 q(a_j+1)>\gamma_j^+-1,&1\le q<\infty,\\
 q(a_j+2)>\gamma_j^+,&0<q<1,
 \end{cases}\qquad j=1,\ldots,n.
\end{equation}
Let $b\in H(\D^n)\cap L^1(\nu_{\aA})$, and fix any positive
integers $m_j>1+\gamma_j^+/p$.
\begin{enumerate}[label=\textup{(\roman*)},leftmargin=*]
\item If $q<p$, put $s=pq/(p-q)$. Then
\[
 \Hb:B_p(\om)\longrightarrow\overline{B_q(\om)}
 \text{ is bounded}\quad\Longleftrightarrow\quad b\in B_s(\om),
\]
and its operator quasi-norm is comparable to $\norm b_{B_s(\om)}$.
Equivalently, the symbol satisfies
\begin{equation}\label{eq:symbolintegral}
 \int_{\D^n}|\T b(\z)|^s\bd{\z}{s\m}v(\z)\,dV(\z)<\infty.
\end{equation}
\item If $p\le q$, then
the same operator is bounded if and only if
\begin{equation}\label{eq:symbolsup}
 \mathcal N_{p,q}(b)
 =\sup_{\z\in\D^n}\bd{\z}{\m}|\T b(\z)|
                    \ow{\z}^{1/q-1/p}<\infty.
\end{equation}
Its operator quasi-norm is comparable to $\mathcal N_{p,q}(b)$.
\end{enumerate}
Constants may depend on the exponents, dimension, weight constants and
the fixed parameters $\aA,\m$, but not on $b$.
\end{theorem}

The derivative order $\m$ is auxiliary: every choice satisfying the
stated lower bounds gives the same symbol class, with equivalent
quantities for each fixed choice. The chosen bound ensures
$pm_j>\gamma_j^+$ for the kernel tests,
$m_j>1+(\gamma_j^+-1)/p$ for synthesis when $p>1$, and
$m_j>(1/p-1/q)\gamma_j^+$ in the range $p\le q$.
The two restrictions in
\eqref{eq:kernelcondition} reflect the two output estimates used in
the proof. For $q>1$ we use a weighted projection bound, at $q=1$
an integral pairing with a growth space, and for $q<1$ the
$q$th-power kernel estimate. The restrictions meet at $q=1$.
For $0<q<1$,
$q(a_j+2)>\gamma_j^+$ has the same form as
$p(\alpha_j+2)>\alpha_{\omega_j}$ in
Harutyunyan--Marinescu \cite[Theorem~2.1 of the preprint]{HM}, after setting $p=q$
and matching the weight exponents. Their theorem treats bounded
measurable symbols with equal input and output exponents. Here we
consider anti-holomorphic symbols and allow the two exponents to differ.
Throughout this paper, \eqref{eq:kernelcondition} is an assumption under
which the characterization holds. We do not claim that it is necessary
for the boundedness of an operator with a fixed analytic symbol.

The diagonal case can be stated in terms of the product Bloch space.
For $S\subseteq\{1,\ldots,n\}$, write
$\partial_S=\prod_{j\in S}\partial_{z_j}$. The product Bloch space
$\Bp$ is defined by the norm
\begin{equation}\label{eq:Blochdef}
 \norm b_{\Bp}=|b(\mathbf0)|+
 \max_{\varnothing\ne S\subseteq\{1,\ldots,n\}}
 \sup_{\z_S\in\D^S}
 \prod_{j\in S}(1-|z_j|^2)
 |\partial_Sb(\z_S,\mathbf0_{S^c})|.
\end{equation}
When $p=q$, the weight in the supremum condition disappears, and
$\mathcal N_{p,p}(b)\asymp\norm b_{\Bp}$ by Lemma~\ref{lem:bloch}.
In particular, for $p=q>1$, the kernel condition can be written as
$p'(a_j+2)-(p'-1)\gamma_j^+>1$.

Section~\ref{sec:model} passes to equivalent Bergman norms using
\cite[Theorem~1]{PR} and establishes the kernel and synthesis estimates.
In Section~\ref{sec:dual}, the output pairing determines a product
space for the Hankel form. Constantin's decomposition
\cite[Corollary~3.2]{Constantin} then allows each atom of this space
to be split into an input factor and a dual factor. This is the
factorization needed for necessity when $q\ge1$.
Section~\ref{sec:smallq} proves the direct output and testing estimates
for $q<1$. Section~\ref{sec:proof} uses these results to prove the
main theorem and its product Bloch corollary.

\section{Kernel estimates, derivative identities and atomic synthesis}\label{sec:model}

For a positive radial product density $w$, let
$A_w^t=H(\D^n)\cap L^t(w\,dV)$.
We first estimate the norms of power kernels and relate $A_v^t$
to $B_t(\om)$. We then establish the synthesis bound that controls
the factors of the atoms and the sums used to test the operator.

\begin{lemma}\label{lem:kernel}
Let $W$ be a positive measurable function on $(0,1]$ such that
\[
 C^{-1}(s/t)^{\theta_+}\le W(s)/W(t)\le C(s/t)^{\theta_-},
 \qquad 0<s\le t\le1,\quad 1<\theta_-\le\theta_+.
\]
Suppose that $w(x)\asymp W(1-x)(1-x)^{-2}$ for $0\le x<1$.
Then $\widehat w(x)=\int_x^1w(u)\,du\asymp W(1-x)/(1-x)$.
If $L>\theta_+$ and $\delta=1-|\lambda|$, then
\begin{equation}\label{eq:kernelestimate}
 \int_\D\frac{w(|z|)}{|1-\overline\lambda z|^L}\dd A(z)
 \asymp W(\delta)\delta^{-L}.
\end{equation}
The corresponding estimates for product densities follow by iteration.
\end{lemma}
\begin{proof}
Put $h=1-x$. Substituting $t=1-u$ and using the comparison for $w$, we have
\[
 \widehat w(x)\asymp\int_x^1\frac{W(1-u)}{(1-u)^2}\,du
 =\int_0^hW(t)t^{-2}\,dt.
\]
Applying the upper and lower scale bounds for $0<t\le h$ gives
\begin{align*}
 \int_0^hW(t)t^{-2}\,dt
 &\le CW(h)h^{-\theta_-}\int_0^h t^{\theta_--2}\,dt
 =\frac{C}{\theta_--1}\frac{W(h)}h,\\
 \int_0^hW(t)t^{-2}\,dt
 &\ge C^{-1}W(h)h^{-\theta_+}\int_0^h t^{\theta_+-2}\,dt
 =\frac{1}{C(\theta_+-1)}\frac{W(h)}h.
\end{align*}
Both power integrals converge because $\theta_+\ge\theta_->1$.
Consequently, $\widehat w(x)\asymp W(h)/h=W(1-x)/(1-x)$,
with comparison constants independent of $x$.

For the upper bound in \eqref{eq:kernelestimate}, apply
Harutyunyan--Lusky \cite[Lemma~2(a)]{HL} with $a=-2$ and $b=L$.
The weight indices in their notation are
$\alpha_W=\theta_+$ and $\beta_W=-\theta_-$, so the required
inequalities follow from $\theta_->1$ and $L>\theta_+$.
To apply their result under the scale assumptions in the statement,
first pass to an equivalent weight with the required representation.
Define
\[
 W_0(t)=\inf_{0<u\le1}W(u)
 \max\left\{(t/u)^{\theta_-},(t/u)^{\theta_+}\right\}.
\]
The scale bounds imply $C^{-1}W(t)\le W_0(t)\le W(t)$.
Moreover, $W_0$ satisfies the same scale bounds with constant one.
Thus $\log W_0(t)$ is locally absolutely continuous as a function
of $\log t$, with logarithmic derivative between $\theta_-$ and
$\theta_+$ almost everywhere. Consequently,
$W_0(t)=W_0(1)\exp(\int_t^1\varepsilon_0(u)\,du/u)$, where
$-\theta_+\le\varepsilon_0\le-\theta_-$.
We may apply the cited lemma to $W_0$ and then replace it by the
comparable weight $W$.
Since $1-r^2\asymp1-r$, the scale bounds also give
$W(1-r^2)\asymp W(1-r)$. It follows that
\[
 \int_\D\frac{w(|z|)}{|1-\overline\lambda z|^L}\dd A(z)
 \lesssim
 \int_\D\frac{W(1-|z|^2)}
 {(1-|z|^2)^2|1-\overline\lambda z|^L}\dd A(z)
 \lesssim W(\delta)\delta^{-L}.
\]

For the reverse inequality, let
$E_\lambda=D(\lambda,\delta/4)\subset\D$.
For $z\in E_\lambda$, we have
\[
 \frac{3\delta}{4}\le1-|z|\le\frac{5\delta}{4},
 \qquad
 \delta\le|1-\overline\lambda z|
 \le1-|\lambda|^2+|\lambda|\,|z-\lambda|
 \le\frac{9\delta}{4}.
\]
The scale bounds imply $W(1-|z|)\asymp W(\delta)$, so that
$w(|z|)\asymp W(\delta)\delta^{-2}$ on $E_\lambda$.
Since $A(E_\lambda)\asymp\delta^2$, we obtain
\begin{align*}
 \int_\D\frac{w(|z|)}{|1-\overline\lambda z|^L}\dd A(z)
 &\ge\int_{E_\lambda}
       \frac{w(|z|)}{|1-\overline\lambda z|^L}\dd A(z)\\
 &\gtrsim W(\delta)\delta^{-L-2}A(E_\lambda)
 \asymp W(\delta)\delta^{-L}.
\end{align*}
This completes the proof of \eqref{eq:kernelestimate}.
For a product weight $w(\z)=\prod_j w_j(z_j)$, Tonelli's theorem
reduces the assertion to the one-variable estimate:
\[
 \int_{\D^n}\prod_j\frac{w_j(z_j)}
 {|1-\overline{\lambda_j}z_j|^{L_j}}\dd V(\z)
 =\prod_j\int_\D\frac{w_j(z_j)}
 {|1-\overline{\lambda_j}z_j|^{L_j}}\dd A(z_j)
 \asymp\prod_j W_j(1-|\lambda_j|)(1-|\lambda_j|)^{-L_j}.
 \qedhere
\]
\end{proof}

The derivative norm formula we need is due to Pel\'aez and de la Rosa.
For a radial weight $w$, write $\widehat w(r)=\int_r^1w(s)\,ds$.
The class $\mathcal D$ consists of the radial weights for which
$\widehat w$ is doubling and reverse doubling in the boundary distance.

\begin{theorem}[Pel\'aez--de la Rosa
  {\cite[Theorem~1, p.~1112]{PR}}]\label{thm:momentLP}
Let $0<t<\infty$ and let $w,\mu\in\mathcal D$ be radial weights
on the unit disc. For $f(z)=\sum_{k\ge0}f_kz^k$, put
$D^\mu f(z)=\sum_{k\ge0}f_kz^k/\mu_{2k+1}$, where
$\mu_x=\int_0^1u^x\mu(u)\,du$. Then
\begin{equation}\label{eq:externalLP}
 \int_\D |f|^tw\,dA
 \asymp\int_\D |D^\mu f|^t\widehat\mu^t w\,dA,
 \qquad f\in H(\D).
\end{equation}
\end{theorem}

The following consequence identifies the Besov spaces with weighted
Bergman spaces and fixes the derivative normalization used below.
\begin{lemma}\label{lem:model}
For every $t>0$,
\begin{equation}\label{eq:model}
 \norm f_{B_t(\om)}\asymp\norm f_{A_v^t}.
\end{equation}
If $a_j>-1$ and each $m_j$ is a positive integer, then
\begin{equation}\label{eq:LP}
 \norm f_{A_v^t}^t\asymp
 \int_{\D^n}|\T f|^t\bd{\z}{t\m} v\dd V.
\end{equation}
The last equivalence also holds with a standard product density
$\bd{\z}{\mathbf c}$, $c_j>-1$, in place of $v$.
\end{lemma}
\begin{proof}
Lemma~\ref{lem:kernel} gives
$\widehat v_j(x)\asymp\omega_j(1-x)/(1-x)$.
The scale estimates imply doubling, and
\[
 \frac{\widehat v_j(1-(1-x)/K)}{\widehat v_j(x)}
 \lesssim K^{1-\gamma_j^-}.
\]
For sufficiently large fixed $K$, the right-hand side is less than
one. Hence $v_j\in\mathcal D$; standard radial densities also belong
to $\mathcal D$.

For \eqref{eq:model}, take $\mu(u)=2$ in
Theorem~\ref{thm:momentLP}. Then $\mu_{2k+1}=1/(k+1)$,
$D^\mu=I+R$, and $\widehat\mu(x)=2(1-x)$.
Applying the theorem successively in each coordinate gives
\eqref{eq:model}, since $1-|z_j|\asymp1-|z_j|^2$.

For \eqref{eq:LP}, use in coordinate $j$ the weight
\[
 \mu_j(u)=\frac{2u^{2a_j+2}(1-u^2)^{m_j-1}}{B(a_j+2,m_j)}.
\]
It belongs to $\mathcal D$, with
$\widehat\mu_j(x)\asymp(1-x)^{m_j}$. The beta integral gives
\[
 (\mu_j)_{2k+1}
 =\frac{B(k+a_j+2,m_j)}{B(a_j+2,m_j)}
 =\frac{(a_j+2)_k}{(a_j+m_j+2)_k},
 \qquad
 \frac1{(\mu_j)_{2k+1}}
 =\prod_{\ell=0}^{m_j-1}\frac{k+a_j+2+\ell}{a_j+2+\ell}.
\]
Thus the product of the coordinate moment derivatives is precisely
$\T$. Another coordinatewise application of
Theorem~\ref{thm:momentLP} proves \eqref{eq:LP}. Tonelli's theorem
justifies the coordinate integrations, also when $t<1$ or either
side is infinite. The same argument applies to standard product densities.
\end{proof}

For a product density $w$ satisfying the hypotheses of
Lemma~\ref{lem:kernel} in each coordinate, the submean inequality
and local comparability give the point evaluation estimate
\begin{equation}\label{eq:point}
 |f(\z)|[w(\z)\bd{\z}{2}]^{1/t}\lesssim\norm f_{A_w^t}.
\end{equation}

Whenever the integral is absolutely convergent, write
$\pair H b=\int H\overline b\,d\nu_{\aA}$.
The next identity transfers derivatives to the symbol while shifting
the weight in the pairing. It follows from the explicit moments of
$\nu_{\aA}$ and is the standard-weight counterpart of the
moment calculation in \cite[Lemma~8 and (2.5) of the preprint]{Wang}.
\begin{lemma}\label{lem:shift}
If $b\in H(\D^n)\cap L^1(\nu_{\aA})$, then
$\T b\in L^1(\nu_{\aA+\m})$ and
\begin{equation}\label{eq:shift}
 \pair H b=\int_{\D^n}H\overline{\T b}\dd\nu_{\aA+\m},
 \qquad H\in\mathcal H_{\mathrm{cl}}.
\end{equation}
Also, for every $\la\in\D^n$,
\begin{equation}\label{eq:kernelpair}
 \pair{\displaystyle\prod_j(1-\overline{\lambda_j}w_j)^{-a_j-m_j-2}}{b}
 =\overline{\T b(\la)}.
\end{equation}
\end{lemma}
\begin{proof}
The standard-weight case of \eqref{eq:LP}, with exponent one, gives $\int|\T b|\,d\nu_{\aA+\m}\lesssim\int|b|\,d\nu_{\aA}$.
In one variable the normalized moment is
\[
 M_a(k)=\int_\D|w|^{2k}\,d\nu_a(w)=\frac{k!}{(a+2)_k}.
\]
The multiplier of $T_{a,m}$ on $w^k$ is
$(a+m+2)_k/(a+2)_k=M_a(k)/M_{a+m}(k)$.
Indeed,
\[
 M_a(k)=2(a+1)\int_0^1x^{2k+1}(1-x^2)^a\,dx
 =(a+1)B(k+1,a+1)=\frac{k!\,\Gamma(a+2)}{\Gamma(k+a+2)}.
\]
For $H\in\mathcal H_{\mathrm{cl}}$, write its coefficients as
$H_{\boldsymbol\beta}$ and those of $b$ as $b_{\boldsymbol\beta}$.
Angular integration now expresses both pairings in terms of the same
moments:
\begin{align*}
 \int H\overline{\T b}\,d\nu_{\aA+\m}
 &=\sum_{\boldsymbol\beta\ge\mathbf0}
 H_{\boldsymbol\beta}\overline{b_{\boldsymbol\beta}}
 \prod_j\left[
 \frac{M_{a_j}(\beta_j)}{M_{a_j+m_j}(\beta_j)}
 M_{a_j+m_j}(\beta_j)\right]\\
 &=\sum_{\boldsymbol\beta\ge\mathbf0}
 H_{\boldsymbol\beta}\overline{b_{\boldsymbol\beta}}
 \prod_jM_{a_j}(\beta_j)
 =\pair H b.
\end{align*}
The geometric decay of the coefficients of $H$, together with the
preceding $L^1$ bounds, ensures absolute convergence of these sums.
For \eqref{eq:kernelpair}, apply \eqref{eq:shift} to the shifted
kernel and use reproduction for $\nu_{\aA+\m}$. Equivalently,
multiply its power-series coefficients by $M_{a_j}(k_j)$ in each
coordinate. The kernel is bounded when its evaluation point is fixed,
so the preceding $L^1$ bounds justify these interchanges.
\end{proof}

Applying the shifted pairing to the kernel gives, for
$F\in\mathcal H_{\mathrm{cl}}$,
\begin{equation}\label{eq:shiftoperator}
 (\Hb F)(\z)=\int_{\D^n}F(\w)\overline{\T b(\w)}
                  K_{\aA}(\w,\z)\,d\nu_{\aA+\m}(\w).
\end{equation}
To see this, take $H(\w)=F(\w)K_{\aA}(\w,\z)$ in
\eqref{eq:shift}. For fixed $\z\in\D^n$ this function belongs
to $\mathcal H_{\mathrm{cl}}$, and the integrals converge because
$\T b\in L^1(\nu_{\aA+\m})$.

For the discrete estimates, use the product Whitney decomposition
of \cite[Proposition~1.1 and the preceding construction, p.~5
of the preprint]{HM}; a finer version is given in
\cite[Example~1, p.~541]{Constantin}. Fix a partition of the disc
into cells whose diameters
are comparable to their distances from the boundary, and take
its product in $\D^n$. Denote the resulting cells by $Q_k$.
Choose enlarged cells $Q_k^*$, still compactly contained in
$\D^n$, with uniformly bounded overlap. For any representatives
$\la_k\in\overline{Q_k}$, put $\delta_{k,j}=1-|\lambda_{k,j}|^2$.
The cells can be chosen so that, uniformly in these representatives,
\begin{gather}
 1-|z_j|^2\asymp\delta_{k,j},\quad
 v(\z)\asymp v(\la_k),\quad
 V(Q_k)\asymp\prod_j\delta_{k,j}^2,
 \qquad \z\in Q_k^*,\label{eq:geometry}\\
 \int_{Q_k}v\,dV\asymp\ow{\la_k},\qquad
 \sum_k\mathbf1_{Q_k^*}\le C.\label{eq:massoverlap}
\end{gather}
Each point of $Q_k$ has a product of small Euclidean discs inside
$Q_k^*$, with coordinate radii comparable to $\delta_{k,j}$.
The submean inequality therefore gives, for every holomorphic
$h$ and $t>0$,
\begin{equation}\label{eq:localmean}
 \sup_{\overline{Q_k}}|h|^t
 \lesssim\prod_j\delta_{k,j}^{-2}\int_{Q_k^*}|h|^t\,dV.
\end{equation}
The same statement holds for anti-holomorphic $h$ by conjugation.
Also, for $\w\in Q_k^*$ and $\z\in\D^n$,
\begin{equation}\label{eq:denomcompare}
 |1-w_j\overline{z_j}|\asymp
 |1-\lambda_{k,j}\overline{z_j}|.
\end{equation}
The notation $Q_k$ enumerates these product cells in a single sequence.
For a locally comparable product density
$w$, put $m_{w,k}=w(\la_k)\bd{\la_k}{2}$, which is comparable
to the $w$-mass of $Q_k$. The next estimate controls sums of
normalized kernels uniformly over the representatives
$\la_k\in\overline{Q_k}$. The fineness of the grid remains fixed.

We use B\'ekoll\'e's weighted projection theorem \cite{Bekolle},
in the form stated in \cite[Theorem~3.1, p.~540]{Constantin}.
The proposition below records its consequence for the present product
weights; its proof only verifies the weight hypothesis.

\begin{proposition}\label{prop:positiveprojection}
Let $1<t<\infty$, and suppose that $w=\prod_jw_j$ satisfies
Lemma~\ref{lem:kernel} in each coordinate, with mass exponents
$1<\theta_{j,-}\le\theta_{j,+}$.
If $\alpha_j>-1$ and $t(\alpha_j+1)>\theta_{j,+}-1$, then
\[
 P_{\boldsymbol\alpha}^{+}h(\z)
 =\int_{\D^n}|h(\boldsymbol\zeta)|
   \prod_j\frac{(1-|\zeta_j|^2)^{\alpha_j}}
   {|1-\overline{\zeta_j}z_j|^{\alpha_j+2}}
   \,dV(\boldsymbol\zeta)
\]
is bounded on $L^t(w)$.
\end{proposition}
\begin{proof}
We verify the one-variable weight condition for B\'ekoll\'e's theorem.
Put $t'=t/(t-1)$ and set $U(z)=w(|z|)/(1-|z|^2)^\alpha$ and
$dA_\alpha=(1-|z|^2)^\alpha dA$. For a boundary box $S$ of scale
$h$, the scale bounds give
\[
 \int_S U\,dA_\alpha\asymp W(h),\qquad
 \int_S U^{-1/(t-1)}\,dA_\alpha
 \asymp h^{t'(\alpha+2)}W(h)^{-1/(t-1)}.
\]
The first integral converges because $\theta_->1$, and the second
because $t(\alpha+1)>\theta_+-1$. Consequently,
\[
 \left(\int_S U\,dA_\alpha\right)
 \left(\int_S U^{-1/(t-1)}\,dA_\alpha\right)^{t-1}
 \lesssim h^{t(\alpha+2)}\asymp A_\alpha(S)^t.
\]
Interior boxes are covered by local comparability. Thus
$U\in B_t(\alpha)$, and the cited theorem gives the one-variable
positive projection bound. Applying these bounds successively in
the coordinates and using Tonelli's theorem proves the assertion.
\end{proof}

\begin{lemma}\label{lem:synthesis}
Let $0<t<\infty$ and let $w$ be a product density as in
Lemma~\ref{lem:kernel}, with coordinate mass exponents $1<\theta_{j,-}\le\theta_{j,+}$.
On such a grid, let
\[
 A_k(\z)=m_{w,k}^{-1/t}
 \prod_j\frac{\delta_{k,j}^{L_j}}
 {|1-\overline{\lambda_{k,j}}z_j|^{L_j}}.
\]
For finite sequences $\mathbf x$,
\begin{equation}\label{eq:synthesis}
 \norm{\sum_k|x_k|A_k}_{L^t(w)}\lesssim\norm{\mathbf x}_{\ell^t}
\end{equation}
provided $L_j>\theta_{j,+}/t$ when $t\le1$, and
$L_j>1+(\theta_{j,+}-1)/t$ when $t>1$.
The same bound holds for analytic kernels and for finitely many
fixed derivative types, normalized in the form
\[
 m_{w,k}^{-1/t}\prod_j
 \frac{\delta_{k,j}^{L_j}z_j^{e_j}}
 {(1-\overline{\lambda_{k,j}}z_j)^{L_j}},
\]
where $\mathbf e$ ranges over a fixed finite set and the kernel
orders satisfy the preceding inequalities.
\end{lemma}
\begin{proof}
For $0<t\le1$, subadditivity and Lemma~\ref{lem:kernel} give
\[
 \int\left(\sum_k|x_k|A_k\right)^tw\,dV
 \le\sum_k|x_k|^t\int A_k^tw\,dV
 \lesssim\sum_k|x_k|^t,
\]
since $tL_j>\theta_{j,+}$ and
$m_{w,k}\asymp\prod_j W_j(1-|\lambda_{k,j}|)$.

Suppose $t>1$ and put $\alpha_j=L_j-2$. Then $\alpha_j>-1$ and
$t(\alpha_j+1)>\theta_{j,+}-1$, so
Proposition~\ref{prop:positiveprojection} applies. Define
\[
 g=\sum_k|x_k|m_{w,k}^{-1/t}\mathbf1_{Q_k}.
\]
The cells are disjoint and their $w$-masses are comparable to
$m_{w,k}$; hence $\|g\|_{L^t(w)}^t\asymp\sum_k|x_k|^t$.
The cell geometry and denominator comparison give, uniformly in
the representatives,
\[
 P_{\boldsymbol\alpha}^+g(\z)
 \asymp\sum_k|x_k|m_{w,k}^{-1/t}
       \prod_j\frac{\delta_{k,j}^{\alpha_j+2}}
       {|1-\overline{\lambda_{k,j}}z_j|^{\alpha_j+2}}
 =\sum_k|x_k|A_k(\z).
\]
The positive projection bound therefore proves \eqref{eq:synthesis}.
Analytic sums are dominated by these positive sums. Factors
$\z^{\mathbf e}$ have modulus at most one, and grouping the fixed
finite set of derivative types only changes the constant.
\end{proof}

Let $\varepsilon_k$ be independent random signs with probabilities
$\mathbb P(\varepsilon_k=1)=\mathbb P(\varepsilon_k=-1)=1/2$.
Then $\mathbb E(\varepsilon_j\varepsilon_k)=\delta_{jk}$, while
Khintchine's inequality \cite[Appendix~A]{Duren} gives, for every
$t>0$ and every finite complex sequence,
\begin{equation}\label{eq:khintchine}
 \left(\mathbb E\left|\sum_k\varepsilon_ka_k\right|^t\right)^{1/t}
 \asymp_t\left(\sum_k|a_k|^2\right)^{1/2}.
\end{equation}
The constants are independent of the number of terms. Applied
pointwise and integrated using Tonelli's theorem, this gives
\begin{equation}\label{eq:khintchinefunctions}
 \mathbb E\int\left|\sum_k\varepsilon_ka_k(\z)\right|^tw(\z)\,dV(\z)
 \asymp_t\int\left(\sum_k|a_k(\z)|^2\right)^{t/2}w(\z)\,dV(\z).
\end{equation}
This is the randomization used in
\cite[Section~3, (3.11)--(3.12) of the preprint]{Wang}.
The factorization below uses only
$\mathbb E(\varepsilon_j\varepsilon_k)=\delta_{jk}$.

\section{Splitting atoms and identifying Hankel forms for \texorpdfstring{$q\ge1$}{q >= 1}}
\label{sec:dual}\label{sec:atomic}
Throughout this section $q\ge1$, and we retain the weight and kernel
assumptions of Theorem~\ref{thm:main}. Pairing the output with a
holomorphic test function turns the operator into a Hankel form.
Since this form depends only on the product of its arguments,
we identify the space of such products and its bounded symbol
functionals. Splitting the atoms supplies the reverse inclusion
in the product-space identity.
\subsection{Recovering the output norm from a pairing}

\begin{proposition}\label{prop:projection}
Let $t>1$ and assume $t(a_j+1)>\gamma_j^+-1$ for all $j$.
The standard Bergman projection $P_{\aA}$ is bounded on $L^t(v)$.
With
$u_t=\bd{\z}{t'\aA}v^{-t'/t},\qquad t'=t/(t-1),$
the pairing $\pair{F}{G}=\int F\overline G\,d\nu_{\aA}$ identifies
$(A_v^t)^*$ with $A_{u_t}^{t'}$ and
$(A_{u_t}^{t'})^*$ with $A_v^t$, with equivalent norms.
\end{proposition}
\begin{proof}
The positive-projection estimate in
Proposition~\ref{prop:positiveprojection}, applied with $w=v$,
gives the boundedness of $P_{\aA}$ on $L^t(v)$.
Write $c_{\aA}=\prod_j(a_j+1)$. The identity
$v^{1/t}u_t^{1/t'}=\bd{\z}{\aA}$ and H\"older's inequality give
\[
 |\pair F G|\le c_{\aA}\norm F_{A_v^t}\norm G_{A_{u_t}^{t'}}.
\]
For surjectivity, extend $\ell\in(A_v^t)^*$ to $L^t(v)$ by
Hahn--Banach and write
$\ell(F)=\int F\overline\varphi v\,dV$, with
$\norm\varphi_{L^{t'}(v)}=\norm\ell$. Set
\[
 \psi=\frac{\varphi v}{c_{\aA}\bd{\z}{\aA}},
 \qquad G=P_{\aA}\psi.
\]
Then $\norm\psi_{L^{t'}(u_t)}=c_{\aA}^{-1}\norm\ell$.
Kernel symmetry transfers the positive-projection bound to
$L^{t'}(u_t)$, and hence $\norm G_{A_{u_t}^{t'}}\lesssim\norm\ell$.
Symmetry and reproduction now give
\[
 \pair F G=\pair{P_{\aA}F}\psi
 =\int F\overline\varphi v\,dV=\ell(F).
\]
The positive-projection bounds justify the interchanges of integrals.
The same argument with the two weights interchanged gives the reverse
duality. Testing against monomials proves uniqueness of the
representing function.
\end{proof}

We now specialize the dual weight to the output space. For $q>1$, put
\begin{equation}\label{eq:Yfinite}
 u(\z)=\prod_j
 \frac{(1-|z_j|^2)^{q'(a_j+2)-2}}
      {\omega_j(1-|z_j|)^{q'-1}}.
\end{equation}
At the endpoint $q=1$, use
\begin{equation}\label{eq:Yendpoint}
 u=\bd{\z}{\aA}/v,\qquad
 A_u^{\infty}=\{G\in H(\D^n):\norm G_{A_u^{\infty}}
      =\sup_{\z\in\D^n}u(\z)|G(\z)|<\infty\}.
\end{equation}
Here $A_u^\infty$ denotes the weighted growth space with norm
$\sup u|G|$; this differs from the $L^\infty(u\,dV)$ norm.

\begin{proposition}\label{prop:outputdual}
Let $1\le q<\infty$ and assume \eqref{eq:kernelcondition}.
For $U\in\overline{A_v^q}$,
\begin{equation}\label{eq:outputdual}
 \norm U_{\overline{A_v^q}}\asymp
 \sup_{\norm G_{A_u^{q\prime}}\le1}
 \left|\int_{\D^n}U G\dd\nu_{\aA}\right|,
 \qquad 1\le q<\infty.
\end{equation}
Dilations are uniformly bounded on $A_u^{q\prime}$. The supremum may be
restricted to functions holomorphic past the closed polydisc, with
a fixed change of constant.
\end{proposition}
\begin{proof}
When $q>1$, the assertion follows from the weighted duality in Proposition~\ref{prop:projection}.
For $q=1$, the upper bound follows directly from
\[
 \left|\int UG\,d\nu_{\aA}\right|
 \le c_{\aA}\int |U||G|u v\,dV
 \le c_{\aA}\norm G_{A_u^{\infty}}\int |U|v\,dV.
\]
For the reverse bound, we construct a representing function for each continuous functional on the
holomorphic space. Write its Hahn--Banach extension as
$\ell(F)=\int F\overline\varphi v\,dV$,
$\norm\varphi_\infty=\norm\ell$. Define
\[
 G(\z)=\int_{\D^n}
     \varphi(\w)\overline{K_{\aA}(\w,\z)}v(\w)\,dV(\w).
\]
The function is holomorphic in $\z$, and its growth is controlled by the kernel estimate in Lemma~\ref{lem:kernel}:
\[
 \int_{\D^n}|K_{\aA}(\w,\z)|v(\w)\dd V(\w)
 \lesssim\frac{v(\z)}{\bd{\z}{\aA}},
\]
Here we use $a_j+2>\gamma_j^+$ at $q=1$. In each coordinate, the right-hand side is
$\omega_j(1-|z_j|)(1-|z_j|)^{-a_j-2}
 \asymp v_j(z_j)(1-|z_j|^2)^{-a_j}$.
Thus $u(\z)|G(\z)|\lesssim\norm\varphi_\infty$.
For $F\in A_v^1$, Fubini and reproduction give
\begin{align*}
 \pair F G
 &=\int_{\D^n}\overline{\varphi(\w)}v(\w)
   \left[\int_{\D^n}F(\z)K_{\aA}(\w,\z)
                     \,d\nu_{\aA}(\z)\right]dV(\w)\\
 &=\int_{\D^n}F(\w)\overline{\varphi(\w)}v(\w)\,dV(\w)
 =\ell(F).
\end{align*}
Absolute integrability for this interchange follows from
\[
 \int |F(\z)|\left[\int |K_{\aA}(\w,\z)|v(\w)dV(\w)\right]
          d\nu_{\aA}(\z)
 \lesssim\int |F(\z)|\frac{v(\z)}{\bd{\z}{\aA}}
                    \bd{\z}{\aA}dV(\z)<\infty.
\]
Thus every continuous functional admits a representative in the stated growth space.
Applying the Banach-space dual norm identity to $\overline U\in A_v^1$ and taking conjugates
proves the reverse bound in \eqref{eq:outputdual}.

For $q=1$, the coordinate factors of $u$ have lower boundary
power exponent $a_j+2-\gamma_j^+>0$. Therefore
$u(\z)/u(\tau\z)\lesssim1$, uniformly for $0<\tau<1$,
and $\norm{G_\tau}_{A_u^{\infty}}\lesssim\norm G_{A_u^{\infty}}$.
Dominated convergence against $|U|v$ shows that these dilations
recover every pairing. For $q>1$, radial dilations are norm contractions.
\end{proof}

\begin{remark}
The endpoint dual depends on the pairing measure.
The space $A_v^1$ is defined using $v\,dV$,
whereas the pairing uses $d\nu_{\aA}$. When these weights agree,
Constantin's duality takes the form
\[
 (A_{\aA}^1)^*\simeq\Bp
 \quad\text{under}\quad
 \langle F,G\rangle_{\aA}=\int F\overline G\,d\nu_{\aA};
\]
see \cite[Proposition~4.3, p.~175]{ConstantinWeak}.
The integral is first defined for $F\in\mathcal H_{\mathrm{cl}}$
and then extended continuously to $A_{\aA}^1$.
Although the cited result uses a common parameter in all coordinates,
its coordinatewise proof applies equally to different $a_j>-1$.
At $q=1$, our
assumption $a_j+2>\gamma_j^+$ excludes that matched-weight
boundary case. When $p=q$, the product test space introduced below is
$A_\eta^r=A_{\aA}^1$ in Section~\ref{sec:atomic}; its dual
will give the product Bloch condition.
\end{remark}

The output pairing turns the Hankel operator into a bilinear form. Define
\begin{equation}\label{eq:bilinear}
 \Lambda_b(F,G)=\int_{\D^n}F G\overline b\dd\nu_{\aA},
 \qquad F,G\in\mathcal H_{\mathrm{cl}},
\end{equation}
and let $\|\Lambda_b\|$ be its best bound by
$\norm F_{A_v^p}\norm G_{A_u^{q\prime}}$ on this class.

\begin{proposition}\label{prop:hankelbilinear}
For $p>0$, $q\ge1$ under \eqref{eq:kernelcondition},
\begin{equation}\label{eq:bilinearop}
 \norm{\Hb}_{A_v^p\to\overline{A_v^q}}\asymp\|\Lambda_b\|,
\end{equation}
with either finite norm determining a bounded extension on the other side.
\end{proposition}
\begin{proof}
For $F,G\in\mathcal H_{\mathrm{cl}}$, put $U=\Hb F$ and
$U_\tau(\z)=U(\tau\z)$. Fubini and reproduction give
\begin{align}
 \int U_\tau(\z)G(\z)\,d\nu_{\aA}(\z)
 &=\int F(\w)\overline{b(\w)}
   \left[\int K_{\aA}(\w,\tau\z)G(\z)
                         \,d\nu_{\aA}(\z)\right]d\nu_{\aA}(\w)\notag\\
 &=\int F(\w)\overline{b(\w)}G(\tau\w)\,d\nu_{\aA}(\w)
 =\Lambda_b(F,G_\tau).\label{eq:dilatedpair}
\end{align}
The same identity holds for $G\in A_u^{q\prime}$. Indeed,
the weight identities and H\"older's inequality give
$A_u^{q\prime}\hookrightarrow L^1(\nu_{\aA})$, with the
growth-space bound used at $q=1$; the dilated kernel is bounded.
For a bounded operator, letting the dilation parameter tend to one in this identity and using the pairing estimate gives
$\|\Lambda_b\|\lesssim\|\Hb\|$.
Conversely, a bound for the bilinear form controls the dilated output: by \eqref{eq:outputdual} and the uniform dilation bound,
\begin{align*}
 \|U_\tau\|_{A_v^q}
 &\lesssim\sup_{\|G\|_{A_u^{q\prime}}\le1}
       |\Lambda_b(F,G_\tau)|\\
 &\le\|\Lambda_b\|\|F\|_{A_v^p}
       \sup_{\|G\|_{A_u^{q\prime}}\le1}\|G_\tau\|_{A_u^{q\prime}}
 \lesssim\|\Lambda_b\|\|F\|_{A_v^p}.
\end{align*}
Fatou's lemma gives the same estimate for $U$. The resulting bound
extends the operator to $A_v^p$ and proves \eqref{eq:bilinearop}.
\end{proof}

\subsection{Decomposing functions and splitting each atom}

H\"older's inequality determines the exponent and weight of the space containing products of the two arguments. Set
\begin{equation}\label{eq:eta}
 d=\frac1p-\frac1q,\qquad
 r=\frac1{1+d},\qquad \eta=\bd{\z}{r\aA}v^{rd}.
\end{equation}
Since $1+d=1/p+1/q'>0$, the exponent $r$ is positive and finite. If $q>1$,
$1/r=1/p+1/q'$ and $\eta=v^{r/p}u^{r/q'}$.
If $q=1$, then $r=p$ and $\eta=vu^p$.
For $q>1$, $p/r>1$, $q'/r>1$ and $r/p+r/q'=1$.
These relations give
\begin{align*}
 v^{r/p}u^{r/q'}
 &=v^{r/p}\left(\bd{\z}{q'\aA}v^{-q'/q}\right)^{r/q'}
 =\bd{\z}{r\aA} v^{r/p-r/q}=\eta,\\
 \int|FG|^r\eta\,dV
 &=\int (|F|^pv)^{r/p}(|G|^{q'}u)^{r/q'}\,dV\\
 &\le\left(\int|F|^pv\,dV\right)^{r/p}
      \left(\int|G|^{q'}u\,dV\right)^{r/q'}.
\end{align*}
For $q=1$, $r=p$ and $\eta=vu^p$, hence
\[
 \int|FG|^p\eta\,dV
 =\int|F|^pv\,(|G|u)^p\,dV
 \le\norm G_{A_u^{\infty}}^p\norm F_{A_v^p}^p.
\]
Taking $r$th roots, or $p$th roots in the endpoint case, gives
\begin{equation}\label{eq:productbound}
 \norm{FG}_{A_\eta^r}\le C\norm F_{A_v^p}\norm G_{A_u^{q\prime}}.
\end{equation}
The normalization of its atoms is determined by the local mass identity
\begin{equation}\label{eq:massidentity}
 (\eta \bd{\z}{2})^{1/r}=\bd{\z}{\aA}\bd{\z}{2}\ow{\z}^d.
\end{equation}
Since $1/r=1+d$ and $v=\ow{\z}/\bd{\z}{2}$, we have
\begin{align*}
 (\eta \bd{\z}{2})^{1/r}
 &=\bd{\z}{\aA}v^d\bd{\z}{2/r}\\
 &=\bd{\z}{\aA}\ow{\z}^d \bd{\z}{2(1/r-d)}\\
 &=\bd{\z}{\aA+2}\ow{\z}^d.
\end{align*}
For $p=q$, $r=1$ and $\eta=\bd{\z}{\aA}$, so
$A_\eta^r=A_{\aA}^1$ with equivalent norms.

\begin{lemma}\label{lem:etaweight}
Each coordinate factor of $\eta$ is integrable, has finite two-sided
boundary power exponents, and has lower density exponent greater
than $-1$. The same statements hold for $u$ when $q>1$.
Consequently, Constantin's polydisc atomic-decomposition theorem applies to these weights,
after replacing them by comparable continuous weights if necessary.
\end{lemma}
\begin{proof}
In coordinate $j$, writing $t=1-|z_j|$, we have
$\eta_j(1-t)\asymp t^{ra_j-2rd}\omega_j(t)^{rd}$.
If $d\ge0$, its lower density exponent plus one is
$\alpha_{\eta,j}^-+1 =r\{a_j+1+d(\gamma_j^--1)\}>0.$
If $d<0$, put $1/s=-d$. Here $p>q\ge1$, $s>q$ and $r=s'$.
In this case the corresponding expression becomes
$\alpha_{\eta,j}^-+1 =r\{a_j+1-(\gamma_j^+-1)/s\}>0$
by \eqref{eq:kernelcondition}. Upper exponents are finite by
\eqref{eq:scale}. For $u$, the lower exponent plus one is
$[q(a_j+1)-(\gamma_j^+-1)]/(q-1)>0$.

If a radial density $w$ has finite exponents
$-1<\alpha_-\le\alpha_+$, choose $P>\max\{1,r\}$ and
$\tau>-1$ so large that $P(\tau+1)>\alpha_++1$.
Then $w/(1-|z|^2)^\tau$ belongs to $B_P(\tau)$.
Indeed, on a boundary box of scale $h$, the two integrals in the
B\'ekoll\'e condition are comparable to
$w(1-h)h^2,\qquad w(1-h)^{-1/(P-1)}h^{P\tau/(P-1)+2},$
respectively. Their product, with the second factor raised to $P-1$,
is $h^{P(\tau+2)}$. The strict inequalities on the exponents justify
both integrations. Finally,
$\widehat w(x)/(1-x)\asymp w(x)$ is a continuous positive
replacement when $w$ is only measurable. This replacement preserves both
the hypotheses and the norm estimates used below.
\end{proof}

We use the scalar-valued polydisc atomic decomposition of Constantin.
Choose a sufficiently fine product grid as in Section~\ref{sec:model};
its fineness remains fixed.

\begin{proposition}[Constantin {\cite[Corollary~3.2]{Constantin}}]\label{prop:atomic}
There are product sampling points $\la_k$ and a fixed finite set of
derivative indices $\mathbf e_k$ such that every $H\in A_\eta^r$
has a norm-convergent expansion
\begin{equation}\label{eq:atomic}
 H=\sum_k c_k\Phi_k,\qquad
 \norm{\mathbf c}_{\ell^r}\lesssim\norm H_{A_\eta^r}.
\end{equation}
Put $\delta_{k,j}=1-|\lambda_{k,j}|^2$ and
$m_{w,k}=w(\la_k)\bd{\la_k}{2}$. After uniformly comparable rescalings
of the coefficients, the atoms may be written
\begin{equation}\label{eq:atomform}
 \Phi_k(\z)=m_{\eta,k}^{-1/r}
 \prod_j\frac{\delta_{k,j}^{M_{k,j}}z_j^{e_{k,j}}}
 {(1-\overline{\lambda_{k,j}}z_j)^{M_{k,j}}},
 \qquad M_{k,j}=\sigma_j+2+e_{k,j}.
\end{equation}
The base orders $\sigma_j$ can be chosen arbitrarily large subject
to fixed lower bounds. The index $k$ includes both a grid index and
a member of the finite derivative set.
\end{proposition}
Here the weight hypotheses follow from Lemma~\ref{lem:etaweight}.
In \cite[Corollary~3.2, p.~543]{Constantin}, take the exponent $r$,
choose $P_j>\max\{1,r\}$ and sufficiently large $\tau_j$, and put
\[
 \sigma_j=\frac{P_j}{r}(\tau_j+2)-2,
 \qquad N_j>2\left(\frac{P_j}{r}-1\right).
\]
The derivative indices range over a fixed finite set determined by
$\mathbf N$. The normalization in \eqref{eq:atomform} follows from
\cite[Lemma~3.1, p.~542]{Constantin}:
\[
 \left\|\overline\partial_{\la}^{\mathbf e}
 K_{\la}^{\boldsymbol\sigma}\right\|_{A_\eta^r}
 \asymp [\eta(\la)\bd{\la}{2}]^{1/r}
 \prod_j(1-|\lambda_j|^2)^{-\sigma_j-2-e_j}.
\]
Increasing $\boldsymbol\tau$ increases the base kernel orders while
leaving the finite derivative set fixed. Comparable continuous
replacements of the weights give the same norms and estimates.

At $q=1$, the second factor lies in a growth space. The following
estimate provides the corresponding synthesis bound.

\begin{lemma}\label{lem:infinitysynthesis}
For $q=1$, normalized analytic kernels
\[
 G_k(\z)=u(\la_k)^{-1}
 \prod_j\frac{\delta_{k,j}^{L_j}}
 {(1-\overline{\lambda_{k,j}}z_j)^{L_j}}
\]
satisfy $\norm{\sum_k y_kG_k}_{A_u^{\infty}}\lesssim\norm{\mathbf y}_{\ell^\infty}$
for finite sums if $L_j>1+a_j+2-\gamma_j^-$.
The assertion also holds for finitely many fixed choices of the orders.
\end{lemma}
\begin{proof}
The coordinate factors of $u$ have boundary exponents between
$\epsilon_{j,-}=a_j+2-\gamma_j^+>0$ and
$\epsilon_{j,+}=a_j+2-\gamma_j^-$.
By comparability within each cell, $\sum_k|G_k(\z)|$ is bounded by a constant
times a product integral whose coordinate integrand is
\[
 \frac{u_j(\zeta)^{-1}(1-|\zeta|)^{L_j-2}}
 {|1-\overline\zeta z_j|^{L_j}}.
\]
Its mass function is $W(t)=t^{L_j}/u_j(1-t)$.
Its lower exponent $L_j-\epsilon_{j,+}$ exceeds one, and its
upper exponent $L_j-\epsilon_{j,-}$ is strictly less than $L_j$.
The exponents thus satisfy the hypotheses of Lemma~\ref{lem:kernel}, which bounds the integral by
$Cu_j(z_j)^{-1}$. Multiplying these bounds over the coordinates proves the estimate.
\end{proof}

We split each atom by distributing its kernel order between two factors, normalized in the input and dual spaces, respectively.

\begin{proposition}\label{prop:split}
The orders in Proposition~\ref{prop:atomic} can be chosen so that
$\Phi_k=F_kG_k$, where $F_k,G_k\in\mathcal H_{\mathrm{cl}}$ and
\begin{equation}\label{eq:factorbounds}
 \norm{\sum_k x_kF_k}_{A_v^p}\lesssim\norm{\mathbf x}_{\ell^p},
 \qquad
 \begin{cases}
 \norm{\sum_k y_kG_k}_{A_u^{q\prime}}\lesssim\norm{\mathbf y}_{\ell^{q'}},&q>1,\\
 \norm{\sum_k y_kG_k}_{A_u^{\infty}}\lesssim\norm{\mathbf y}_{\ell^\infty},&q=1.
 \end{cases}
\end{equation}
In particular both individual factor norms are uniformly bounded.
\end{proposition}
\begin{proof}
Choose $L_j$ large enough for Lemma~\ref{lem:synthesis} with
weight $v$ and exponent $p$. Then increase $\sigma_j$ until every
$M_{k,j}-L_j$ satisfies the synthesis condition for the second factor.
Define
\[
 F_k(\z)=m_{v,k}^{-1/p}
 \prod_j\frac{\delta_{k,j}^{L_j}z_j^{e_{k,j}}}
 {(1-\overline{\lambda_{k,j}}z_j)^{L_j}}.
\]
For $q>1$, take
\[
 G_k(\z)=m_{u,k}^{-1/q'}
 \prod_j\frac{\delta_{k,j}^{M_{k,j}-L_j}}
 {(1-\overline{\lambda_{k,j}}z_j)^{M_{k,j}-L_j}}.
\]
The normalizations agree because $m_{\eta,k}^{1/r}=m_{v,k}^{1/p}m_{u,k}^{1/q'}$.
Indeed,
\begin{align*}
 m_{v,k}^{1/p}m_{u,k}^{1/q'}
 &=v(\la_k)^{1/p}u(\la_k)^{1/q'}
                    \bd{\la_k}{2(1/p+1/q')}\\
 &=\eta(\la_k)^{1/r}\bd{\la_k}{2(1/r)}
  =m_{\eta,k}^{1/r}.
\end{align*}
For each coordinate, the product of the two kernel factors is
\[
 \frac{\delta_{k,j}^{L_j}z_j^{e_{k,j}}}
      {(1-\overline{\lambda_{k,j}}z_j)^{L_j}}
 \frac{\delta_{k,j}^{M_{k,j}-L_j}}
      {(1-\overline{\lambda_{k,j}}z_j)^{M_{k,j}-L_j}}
 =\frac{\delta_{k,j}^{M_{k,j}}z_j^{e_{k,j}}}
       {(1-\overline{\lambda_{k,j}}z_j)^{M_{k,j}}}.
\]
Hence $F_kG_k=\Phi_k$.
The local mass normalization can also be seen directly from
\begin{align*}
 \norm{F_k}_{A_v^p}^p
 &\le m_{v,k}^{-1}\prod_j\left[
       \delta_{k,j}^{pL_j}\int_\D
       \frac{v_j(z_j)}{|1-\overline{\lambda_{k,j}}z_j|^{pL_j}}\,dA(z_j)
                                    \right]\\
 &\lesssim m_{v,k}^{-1}\prod_j\omega_j(1-|\lambda_{k,j}|)
 =1.
\end{align*}
For $G_k$, the same calculation with $u,q'$ and $M_{k,j}-L_j$
gives $\norm{G_k}_{A_{u}^{q'}}\lesssim1$.
The sequence estimates \eqref{eq:factorbounds} follow from
Lemma~\ref{lem:synthesis}, applied with the respective weights and exponents to the majorants
\[
 \left|\sum_kx_kF_k\right|\le\sum_k|x_k||F_k|,
 \qquad
 \left|\sum_ky_kG_k\right|\le\sum_k|y_k||G_k|.
\]
For $q=1$, replace $m_{u,k}^{-1/q'}$ by $u(\la_k)^{-1}$.
Now the corresponding normalization is
$m_{\eta,k}^{1/p} =[v(\la_k)u(\la_k)^p\bd{\la_k}{2}]^{1/p} =m_{v,k}^{1/p}u(\la_k).$
Lemma~\ref{lem:infinitysynthesis} gives
\[
 u(\z)\left|\sum_ky_kG_k(\z)\right|
 \le\norm{\mathbf y}_{\ell^\infty}
       u(\z)\sum_k|G_k(\z)|
 \lesssim\norm{\mathbf y}_{\ell^\infty}.
\]
Grouping the terms by derivative type introduces only a fixed factor in the constants, since the number of types is finite.
\end{proof}

To pass from individual atoms to general functions, we measure
the cost of summing these products.
Define the weak product with exponent $\min\{1,r\}$ to consist
of functions with representations
$H=\sum_\ell F_\ell G_\ell$, converging in $A_\eta^r$, and finite
quasi-norm
\begin{equation}\label{eq:weaknorm}
 \norm H_{A_v^p\odot A_u^{q\prime}}
 =\inf\left[\sum_\ell
  (\norm{F_\ell}_{A_v^p}\norm{G_\ell}_{A_u^{q\prime}})^{\min\{1,r\}}
             \right]^{1/\min\{1,r\}}.
\end{equation}
The infimum is over countable product representations, including
finite ones. By \eqref{eq:productbound}, finite cost guarantees
convergence in $A_\eta^r$. Completeness follows by concatenating
representations of a subsequence of successive differences with
summable $\min\{1,r\}$-costs.
For $r\ge1$, this definition agrees with the usual weak product.

For related weak factorizations, see
\cite[Theorem~4.1]{ConstantinWeak} for standard weights on the
polydisc and \cite[Theorem~2 of the preprint]{Wang} for doubling
weights on the disc. We use the preceding polydisc decomposition
to prove the factorization for the present product weights,
including the growth-space factor at $q=1$.

\begin{theorem}\label{thm:factor}
Let $0<p<\infty$ and $1\le q<\infty$. Under the weight and kernel
assumptions of Theorem~\ref{thm:main},
\begin{equation}\label{eq:weakfactor}
 A_\eta^r=A_v^p\odot A_u^{q\prime}
\end{equation}
with equivalent quasi-norms. For every finite atom sum in
\eqref{eq:atomic}, its weak-product quasi-norm is at most
$C\norm{\mathbf c}_{\ell^r}$, with factors in
$\mathcal H_{\mathrm{cl}}$.
\end{theorem}
\begin{proof}
The inclusion of the weak product in the weighted Bergman space follows from H\"older's inequality and the triangle or $r$-triangle inequality, according to
$r\ge1$ or $r<1$:
$\norm H_{A_\eta^r}\lesssim\norm H_{A_v^p\odot A_u^{q\prime}}$.
For the reverse inclusion, write $H=\sum_kc_k\Phi_k$ as in
Proposition~\ref{prop:atomic} and use $\Phi_k=F_kG_k$ from
Proposition~\ref{prop:split}. When $r\le1$, the coefficient norm
controls the sum directly. When $r>1$, we retain the $\ell^r$
bound by averaging finite products.

If $r\le1$, the representation
$\sum_kc_k\Phi_k=\sum_k(c_kF_k)G_k$ has $r$-cost at most
$C\norm{\mathbf c}_{\ell^r}$, because
\[
 \left[\sum_k\big(\norm{c_kF_k}_{A_v^p}
                         \norm{G_k}_{A_u^{q\prime}}\big)^r\right]^{1/r}
 \le C\left(\sum_k|c_k|^r\right)^{1/r}.
\]
This proves the assertion for finite sums in this range.
For scalar pairings, we may also use
\begin{equation}\label{eq:ellsum}
 \sum_k|c_k|\le\left(\sum_k|c_k|^r\right)^{1/r}.
\end{equation}
Indeed, divide by $R=\norm{\mathbf c}_{\ell^r}>0$. Then
$|c_k|/R\le1$, so
$\sum_k|c_k|/R\le\sum_k(|c_k|/R)^r=1$.

If $r>1$ and $q>1$, then $p>q>1$.
For a finite coefficient vector set
$\alpha_k=e^{i\arg c_k}|c_k|^{r/p},\qquad \beta_k=|c_k|^{r/q'},$
with both zero at zero coefficients. Let $\varepsilon_k$ be
independent signs of equal probabilities. Since
$\mathbb E(\varepsilon_j\varepsilon_k)=\delta_{jk}$,
\begin{equation}\label{eq:randomidentity}
 \sum_k c_k\Phi_k
 =\mathbb E\left[
   \left(\sum_k\varepsilon_k\alpha_kF_k\right)
   \left(\sum_k\varepsilon_k\beta_kG_k\right)\right].
\end{equation}
The expectation is a finite average of products in which the mixed terms cancel. Expanding it gives
\begin{align*}
 \mathbb E\left[\left(\sum_j\varepsilon_j\alpha_jF_j\right)
                  \left(\sum_k\varepsilon_k\beta_kG_k\right)\right]
 &=\sum_{j,k}\alpha_j\beta_kF_jG_k
                        \mathbb E(\varepsilon_j\varepsilon_k)\\
 &=\sum_k\alpha_k\beta_kF_kG_k
 =\sum_k e^{i\arg c_k}|c_k|^{r/p+r/q'}\Phi_k
 =\sum_kc_k\Phi_k.
\end{align*}
The final equality uses $r/p+r/q'=1$.
For each choice of signs, the synthesis estimates \eqref{eq:factorbounds} bound the
product of the factor norms by
\[
 C\left(\sum_k|c_k|^r\right)^{1/p+1/q'}
 =C\norm{\mathbf c}_{\ell^r}.
\]
Averaging preserves this bound for the weak-product cost.
More explicitly, writing $N$ for the number of coefficients and
$F_{\boldsymbol\varepsilon}=\sum_k\varepsilon_k\alpha_kF_k$,
$G_{\boldsymbol\varepsilon}=\sum_k\varepsilon_k\beta_kG_k$, we have
\begin{align*}
 \norm{\sum_kc_k\Phi_k}_{A_v^p\odot A_u^{q\prime}}
 &\le2^{-N}\sum_{\boldsymbol\varepsilon\in\{-1,1\}^N}
         \norm{F_{\boldsymbol\varepsilon}}_{A_v^p}
         \norm{G_{\boldsymbol\varepsilon}}_{A_u^{q\prime}}\\
 &\le2^{-N}\sum_{\boldsymbol\varepsilon}
       C\left(\sum_k|c_k|^r\right)^{1/p}
        \left(\sum_k|c_k|^r\right)^{1/q'}
 =C\norm{\mathbf c}_{\ell^r}.
\end{align*}
The constant is independent of $N$, so the same estimate applies
to every finite block of the atomic expansion.

If $q=1$ and $r=p>1$, use instead
\[
 \sum_k c_k\Phi_k
 =\mathbb E\left[
   \left(\sum_k\varepsilon_k c_kF_k\right)
   \left(\sum_k\varepsilon_k G_k\right)\right].
\]
The first norm is controlled by $\norm{\mathbf c}_{\ell^p}$,
and the second is uniformly bounded by the $\ell^\infty$ estimate.
Thus the required weak-product bound holds for every finite atom sum.

For $r\le1$, the original atomic series already has finite
$r$-cost. For $r>1$, choose increasing finite truncations whose
coefficient tail norms are at most $2^{-j}\norm{\mathbf c}_{\ell^r}$.
Partition the expansion into the resulting finite blocks and apply
the finite random representation to each block. The sum of all
block costs is at most $C\norm{\mathbf c}_{\ell^r}$.
By the first embedding and the atomic expansion, the concatenated product series converges to $H$ in $A_\eta^r$.
This representation proves the reverse estimate.
\end{proof}

\begin{remark}
If $p\le q$, then $r\le1$, and splitting the atoms individually
suffices. If $p>q$, then $r>1$, and direct summation would produce
$\sum|c_k|$. For example, $c_k=N^{-1/r}$ for $1\le k\le N$
gives $\norm{\mathbf c}_{\ell^r}=1$ but
$\sum|c_k|=N^{1-1/r}$. The sign average cancels mixed products
and keeps the estimate in terms of the $\ell^r$ norm. The only property of the signs used here is the finite identity
\[
 2^{-N}\sum_{\boldsymbol\varepsilon\in\{-1,1\}^N}
       \varepsilon_j\varepsilon_k=\delta_{jk};
\]
thus expectations can be replaced throughout by finite sign averages.
Khintchine's inequality is not needed for this proof.
\end{remark}

\subsection{Identifying the symbol functional}
We next characterize the symbols that induce bounded functionals on the weak product.
For $H\in\mathcal H_{\mathrm{cl}}$, set $L_b(H)=\pair H b$ and
\begin{equation}\label{eq:pairnorm}
 \|L_b\|=\sup_{\substack{H\in\mathcal H_{\mathrm{cl}}\\
                         \norm H_{A_\eta^r}\le1}}|\pair H b|.
\end{equation}
When this norm is finite, the functional extends uniquely and continuously to $A_\eta^r$.
We retain the same notation for its continuous pairing on $A_\eta^r$;
whenever the shifted integral is absolutely convergent, it represents that pairing.

\begin{proposition}\label{prop:finitepair}
If $p>q\ge1$ and $s=pq/(p-q)$, then
\begin{equation}\label{eq:finitepair}
 \|L_b\|\asymp\norm b_{A_v^s},
\end{equation}
with finiteness on either side equivalent to finiteness on the other.
\end{proposition}
\begin{proof}
Here $r=s'>1$ and $\eta=\bd{\z}{s'\aA}v^{-s'/s}$.
Since $s>q$, condition \eqref{eq:kernelcondition} implies
$s(a_j+1)>\gamma_j^+-1$.
The weighted duality in Proposition~\ref{prop:projection}, applied with $t=s$, therefore identifies
$(A_\eta^{s'})^*$ with $A_v^s$ under our pairing.
If $\|L_b\|<\infty$, let $\widetilde b\in A_v^s$
represent its extension, with $\norm{\widetilde b}_{A_v^s}
\lesssim\|L_b\|$. Testing $H=\z^{\boldsymbol\alpha}$ yields
\[
 \overline{b_{\boldsymbol\alpha}}\prod_jM_{a_j}(\alpha_j)
 =L_b(\z^{\boldsymbol\alpha})
 =\overline{\widetilde b_{\boldsymbol\alpha}}
                         \prod_jM_{a_j}(\alpha_j).
\]
Since all the moments are positive, the coefficients agree:
$b_{\boldsymbol\alpha}=\widetilde b_{\boldsymbol\alpha}$ for every
$\boldsymbol\alpha$. Thus $b=\widetilde b$ and
$\norm b_{A_v^s}\lesssim\|L_b\|$.
The reverse inequality follows directly from H\"older's inequality and the weight identity
\[
 \eta^{1/r}v^{1/s}
 =\bd{\z}{\aA}v^{-1/s}v^{1/s}
 =\bd{\z}{\aA},\qquad \frac1r+\frac1s=1.
\]
Therefore, for every $H\in A_\eta^r$,
\begin{align*}
 |L_b(H)|
 &\le c_{\aA}\int(|H|\eta^{1/r})(|b|v^{1/s})\,dV\\
 &\le c_{\aA}\left(\int|H|^r\eta\,dV\right)^{1/r}
                  \left(\int|b|^sv\,dV\right)^{1/s}.
\end{align*}
Take the supremum over $\norm H_{A_\eta^r}\le1$.
\end{proof}

For $p\le q$, the product exponent satisfies $r\le1$.
We need the weighted $L^1$ embedding familiar from the one-variable
result of Blasco \cite[Lemma~2.1]{Blasco}. For the present product
weights it follows directly from \eqref{eq:point}.

\begin{lemma}\label{lem:embedding}
Suppose $p\le q$, and retain $d,r,\eta$ from \eqref{eq:eta}.
Then $0<r\le1$ and
\begin{equation}\label{eq:embedding}
 \int_{\D^n}|H(\z)|\bd{\z}{\aA}\ow{\z}^d\dd V(\z)
 \lesssim\norm H_{A_\eta^r}.
\end{equation}
\end{lemma}
\begin{proof}
Put $\chi=\bd{\z}{\aA}\ow{\z}^d$. By
\eqref{eq:massidentity},
$\eta=\chi^r\bd{\z}{2r-2}$, and \eqref{eq:point} gives
\[
 |H(\z)|\chi(\z)\bd{\z}{2}
 =|H(\z)|[\eta(\z)\bd{\z}{2}]^{1/r}
 \lesssim\norm H_{A_\eta^r}.
\]
If $r<1$, it follows that
\begin{align*}
 \int|H|\chi\,dV
 &=\int |H|^r\eta
       \bigl(|H|\chi\bd{\z}{2}\bigr)^{1-r}\,dV\\
 &\lesssim\norm H_{A_\eta^r}^{1-r}
       \int|H|^r\eta\,dV
 =\norm H_{A_\eta^r}.
\end{align*}
For $r=1$, the assertion is the identity $\chi=\eta$.
\end{proof}

\begin{proposition}\label{prop:suppair}
If $p\le q$, and $m_j>1+\gamma_j^+/p$, then
\begin{equation}\label{eq:suppair}
 \|L_b\|\asymp\mathcal N_{p,q}(b).
\end{equation}

\end{proposition}
\begin{proof}
Let $d=1/p-1/q\ge0$. To recover the growth of the symbol from a bounded functional, assume $\|L_b\|<\infty$ and test it on
$H_{\la}(\w)=\prod_j(1-\overline{\lambda_j}w_j)^{-a_j-m_j-2}.$
These kernels belong to $\mathcal H_{\mathrm{cl}}$.
The upper mass exponent of $\eta_j$ is
$r(a_j+2+d\gamma_j^+)$, and
$m_j>d\gamma_j^+$ by the chosen bound. Lemma~\ref{lem:kernel}
therefore gives
\[
 \norm{H_{\la}}_{A_\eta^r}
 \lesssim(\eta(\la)\bd{\la}{2})^{1/r}
      \prod_j (1-|\lambda_j|^2)^{-a_j-m_j-2}
 =\ow{\la}^d\bd{\la}{-\m}.
\]
The last equality follows by substituting the local mass identity \eqref{eq:massidentity}:
\[
 [\eta(\la)\bd{\la}{2}]^{1/r}
 =\prod_j (1-|\lambda_j|^2)^{a_j+2}
                  \omega_j(1-|\lambda_j|)^d.
\]
Multiplying by the displayed kernel powers cancels $a_j+2$,
leaving $(1-|\lambda_j|^2)^{-m_j}$.
The kernel pairing \eqref{eq:kernelpair} then recovers the shifted derivative of the symbol at each $\la$:
\begin{align*}
 |\T b(\la)|
 &=|L_b(H_{\la})|
 \le\|L_b\|\norm{H_{\la}}_{A_\eta^r}\\
 &\lesssim\|L_b\|\ow{\la}^{d}\bd{\la}{-\m},\\
 \bd{\la}{\m}|\T b(\la)|\ow{\la}^{-d}
 &\lesssim\|L_b\|.
\end{align*}
Taking the supremum proves $\mathcal N_{p,q}(b)\lesssim\|L_b\|$.

Conversely, suppose that $\mathcal N_{p,q}(b)<\infty$.
For $H\in\mathcal H_{\mathrm{cl}}$, the shifted identity gives
\begin{align*}
 |L_b(H)|
 &\lesssim\int |H|\,|\T b|\bd{\z}{\m}\bd{\z}{\aA}\dd V\\
 &\le\mathcal N_{p,q}(b)
       \int |H|\bd{\z}{\aA}\ow{\z}^d\dd V
 \lesssim\mathcal N_{p,q}(b)\norm H_{A_\eta^r},
\end{align*}
Lemma~\ref{lem:embedding} gives the last inequality. The same
bound proves absolute convergence of the shifted integral for
every $H\in A_\eta^r$, and this integral represents the continuous
extension of $L_b$. Taking the supremum gives the reverse estimate.
\end{proof}

\section{Estimating outputs and testing atoms for \texorpdfstring{$0<q<1$}{0 < q < 1}}
\label{sec:smallq}
For $0<q<1$, we prove sufficiency by estimating the $q$th power
of the output integral. Necessity follows from normalized kernel
tests: individual tests give the supremum condition, while randomized
sums recover the integral condition when $q<p$. On the product cells
from Section~\ref{sec:model}, choose arbitrary
$\la_k\in\overline{Q_k}$ and set
\begin{equation}\label{eq:atoms}
 F_k(\z)=\ow{\la_k}^{-1/p}
 \prod_j\frac{\delta_{k,j}^{m_j}}
 {(1-\overline{\lambda_{k,j}}z_j)^{m_j}}.
\end{equation}
Applying Lemma~\ref{lem:synthesis} with $w=v$, $t=p$ and $L_j=m_j$,
we have $\|\sum_kc_kF_k\|_{A_v^p}\lesssim\|\mathbf c\|_{\ell^p}$;
the required bounds on the kernel powers follow from
$m_j>1+\gamma_j^+/p$. The estimate is uniform in the choice of
representatives and in the signs of the coefficients.

Randomized analytic tests are a standard tool in Bergman-space
embedding problems; see \cite{LueckingKhinchine}. Here we use the
Hankel-form version in \cite[Section~3, (3.11)--(3.12) of the
preprint]{Wang}. We retain the local submean argument below to pass
from the square function to point evaluations on the product cells.

\begin{lemma}\label{lem:randomtest}
If $0<q<1$ and $\Hb:A_v^p\to\overline{A_v^q}$ is bounded with quasi-norm $M$,
then every finite coefficient vector satisfies
\begin{equation}\label{eq:randomtest}
 \sum_k|c_k|^q\ow{\la_k}
                |(\Hb F_k)(\la_k)|^q
 \lesssim M^q\norm{\mathbf c}_{\ell^p}^q.
\end{equation}
\end{lemma}
\begin{proof}
For each choice of signs, the synthesis estimate controls the input
norm, so boundedness of the operator implies
\[
 \int\left|\sum_k\varepsilon_kc_k(\Hb F_k)(\z)\right|^q
       v(\z)\,dV(\z)
 \lesssim M^q\norm{\mathbf c}_{\ell^p}^q.
\]
Interchanging the finite sign average with the integral and applying
\eqref{eq:khintchine} pointwise, we obtain
\begin{equation}\label{eq:square}
 \int\left(\sum_k|c_k(\Hb F_k)(\z)|^2\right)^{q/2}v(\z)\,dV(\z)
 \lesssim M^q\norm{\mathbf c}_{\ell^p}^q.
\end{equation}
To pass from this integral bound to the values at the sampling points,
we use the local submean estimate at $\la_k$:
\[
 \ow{\la_k}|(\Hb F_k)(\la_k)|^q
 \lesssim\int_{Q_k^*}|\Hb F_k|^qv\,dV.
\]
At each point, at most a fixed number $C_0$ of the enlarged cells
overlap. We can therefore compare the sum over these cells with the
square function. H\"older's inequality, applied to the active indices
with exponents $2/q$ and $2/(2-q)$, gives
\[
 \sum_k|c_k(\Hb F_k)(\z)|^q\mathbf1_{Q_k^*}(\z)
 \le C_0^{1-q/2}
       \left(\sum_k|c_k(\Hb F_k)(\z)|^2\right)^{q/2}.
\]
After multiplication by $v$ and integration, \eqref{eq:square}
proves \eqref{eq:randomtest}.
\end{proof}

For sufficiency, we estimate the shifted symbol integral on the
product cells. This uses the cell decomposition and $q$th-power
subadditivity from \cite[proof of Theorem~2.1,
pp.~6--7 of the preprint]{HM}. The calculation below records
the kernel exponent required for the present symbol condition.

\begin{lemma}\label{lem:outputestimate}
Assume $0<q<1$ and \eqref{eq:kernelcondition}. For $F\in A_v^p$
such that the right-hand side below is finite, the shifted integral
\eqref{eq:shiftoperator} is absolutely convergent and satisfies
\begin{equation}\label{eq:outputestimate}
 \norm{\Hb F}_{A_v^q}^q
 \lesssim\int_{\D^n}|F(\w)|^q|\T b(\w)|^q
                  \bd{\w}{q\m}v(\w)\,dV(\w).
\end{equation}

\end{lemma}
\begin{proof}
Put $h=F\T b$, which is holomorphic, and consider the positive integral
\[
 S(\z)=\int_{\D^n}|h(\w)|\,|K_{\aA}(\w,\z)|
                         \,d\nu_{\aA+\m}(\w),
\]
initially allowing the value $+\infty$. We apply the estimates
below first to finite unions of cells and then use monotone
convergence. Local comparability of the kernel and the boundary
factors bounds the integral over each cell $Q_k$ by
\[
 C\sup_{Q_k}|h|\prod_j
    \frac{\delta_{k,j}^{a_j+m_j+2}}
         {|1-\lambda_{k,j}\overline{z_j}|^{a_j+2}}.
\]
The inequality $(\sum x_k)^q\le\sum x_k^q$ now allows us to estimate
the cells separately. Applying \eqref{eq:localmean} to $h$, we obtain
\begin{align*}
 S(\z)^q
 &\lesssim\sum_k\prod_j
     \frac{\delta_{k,j}^{q(a_j+m_j+2)-2}}
          {|1-\lambda_{k,j}\overline{z_j}|^{q(a_j+2)}}
                      \int_{Q_k^*}|h(\w)|^q\,dV(\w)\\
 &\lesssim\int_{\D^n}|h(\w)|^q
       \prod_j\frac{(1-|w_j|^2)^{q(a_j+m_j+2)-2}}
                    {|1-w_j\overline{z_j}|^{q(a_j+2)}}\,dV(\w).
\end{align*}
In the last line, local comparability replaces the sampling points
by integration variables, while bounded overlap controls the sum over
the enlarged cells.

Integrating against $v(\z)$ and applying Tonelli's theorem, we use
the assumption $q(a_j+2)>\gamma_j^+$ to estimate the remaining
kernel integral by Lemma~\ref{lem:kernel}:
\[
 \int_\D\frac{v_j(|z|)}{|1-w_j\overline z|^{q(a_j+2)}}\,dA(z)
 \lesssim\omega_j(1-|w_j|)(1-|w_j|^2)^{-q(a_j+2)}.
\]
Substitution into the preceding estimate cancels the kernel powers
and leaves the desired weight:
\begin{align*}
 \norm S_{L^q(v)}^q
 &\lesssim\int_{\D^n}|F\T b|^q
  \prod_j\big[(1-|w_j|^2)^{q(a_j+m_j+2)-2}
     \omega_j(1-|w_j|)(1-|w_j|^2)^{-q(a_j+2)}\big]\,dV\\
 &=\int_{\D^n}|F\T b|^q\bd{\w}{q\m}v(\w)\,dV.
\end{align*}
The right-hand side is finite by assumption, so $S(\z_0)<\infty$
for some $\z_0\in\D^n$. At this point,
\[
 |K_{\aA}(\w,\z_0)|
 \ge\prod_j(1+|z_{0,j}|)^{-a_j-2}>0,
 \qquad \w\in\D^n.
\]
It follows that $\int|F\T b|\,d\nu_{\aA+\m}<\infty$.
The kernel is uniformly bounded on each compact set of output
points, so the shifted integral converges locally uniformly and
defines an anti-holomorphic function. Since $|\Hb F|\le S$,
the preceding estimate proves \eqref{eq:outputestimate}.

The calculation also explains the agreement with the kernel
restriction in \cite[Theorem~2.1 of the preprint]{HM} when
$p=q<1$. Their derivative estimate has kernel exponent
$p(\alpha_j+3)$ and weight factor $(1-|z_j|^2)^{p-2}$;
the difference is $p(\alpha_j+3)-(p-2)-2=p(\alpha_j+2)$.
Thus the same $q$th-power kernel estimate accounts for the parameter
restriction in both arguments.

\end{proof}

\section{Proving the boundedness criterion and the diagonal case}
\label{sec:proof}
\begin{proof}[Proof of Theorem~\ref{thm:main} for $q\ge1$]
We use the equivalent Bergman norms from Lemma~\ref{lem:model}
and retain $u,\eta,r$ from Section~\ref{sec:dual}. It is enough
to characterize boundedness of $\Hb:A_v^p\to\overline{A_v^q}$.
For sufficiency we estimate the form on products $FG$; for
necessity we use the atomic factorization to control its action
on the whole product space.

\emph{Sufficiency.}
Let $F\in A_v^p$ and $G\in A_u^{q\prime}$.
We estimate the Hankel form in each range of exponents.
Suppose first that $p>q$ and $b\in B_s(\om)$, where
$s=pq/(p-q)$. Then $b\in A_v^s$, $r=s'$, and
$\eta^{1/r}v^{1/s}=\bd{\z}{\aA}$. This identity separates the
weight in the pairing into the two factors needed for H\"older's inequality:
\begin{align*}
 |\Lambda_b(F,G)|
 &=\left|\int_{\D^n}FG\overline b\,d\nu_{\aA}\right|\\
 &\le c_{\aA}\int_{\D^n}
      (|FG|\eta^{1/r})(|b|v^{1/s})\,dV\\
 &\le c_{\aA}\norm{FG}_{A_\eta^r}\norm b_{A_v^s}
 \lesssim\norm b_{A_v^s}\norm F_{A_v^p}\norm G_{A_u^{q\prime}},
\end{align*}
where the last estimate is the product inequality
\eqref{eq:productbound}. The same argument covers $q=1$,
when the second factor carries the weighted supremum norm.

Now suppose that $p\le q$ and $\mathcal N_{p,q}(b)<\infty$.
Here $d=1/p-1/q\ge0$ and $r\le1$. For general $F$ and $G$,
define the form by the shifted integral below, which agrees with
\eqref{eq:bilinear} on $\mathcal H_{\mathrm{cl}}$ by
Lemma~\ref{lem:shift}.
By \eqref{eq:embedding} and the bound
$\bd{\z}{\m}|\T b(\z)|\le\mathcal N_{p,q}(b)\ow{\z}^d$, we obtain
\begin{align*}
 |\Lambda_b(F,G)|
 &=\left|\int_{\D^n}FG\overline{\T b}\,d\nu_{\aA+\m}\right|\\
 &\lesssim\mathcal N_{p,q}(b)
       \int_{\D^n}|FG|\bd{\z}{\aA}\ow{\z}^d\,dV\\
 &\lesssim\mathcal N_{p,q}(b)\norm{FG}_{A_\eta^r}\\
 &\lesssim\mathcal N_{p,q}(b)
       \norm F_{A_v^p}\norm G_{A_u^{q\prime}}.
\end{align*}
The estimate proves absolute convergence of the shifted integral
and boundedness of the extended form.

In both ranges, Proposition~\ref{prop:hankelbilinear} converts the
form bound into the required operator estimate:
\[
 \norm{\Hb F}_{A_v^q}\lesssim\norm F_{A_v^p}
 \begin{cases}
 \norm b_{A_v^s},&p>q,\\
 \mathcal N_{p,q}(b),&p\le q.
 \end{cases}
\]
This holds for every $F\in A_v^p$, and Lemma~\ref{lem:model}
transfers the estimate to the Besov spaces.

\emph{Necessity.}
Assume that $\Hb:A_v^p\to\overline{A_v^q}$ is bounded, and write
$M=\norm{\Hb}_{A_v^p\to\overline{A_v^q}}$.
For $F\in A_v^p$ and $G\in A_u^{q\prime}$, define the bounded
form by pairing the output with $G$. This agrees with the original
Hankel form on $\mathcal H_{\mathrm{cl}}$, while the pairing estimate
in Section~\ref{sec:dual} gives
\begin{align*}
 |\Lambda_b(F,G)|
 &=\left|\int_{\D^n}(\Hb F)(\z)G(\z)\,d\nu_{\aA}(\z)\right|\\
 &\lesssim\norm{\Hb F}_{A_v^q}\norm G_{A_u^{q\prime}}
 \le M\norm F_{A_v^p}\norm G_{A_u^{q\prime}}.
\end{align*}
The remaining step is to pass from products $FG$ to arbitrary
$H\in A_\eta^r$. Since continuity of $L_b$ on this space is
what we need to prove, first take $b_\tau(\z)=b(\tau\z)$,
$0<\tau<1$. The functional $L_{b_\tau}$ is already continuous;
we will obtain a bound independent of $\tau$ from the atomic factors.
Radial dilations are self-adjoint for the $d\nu_{\aA}$ pairing,
and $(FG)_\tau=F_\tau G_\tau$. Hence, for
$F,G\in\mathcal H_{\mathrm{cl}}$,
$\Lambda_{b_\tau}(F,G) =\pair{FG}{b_\tau} =\pair{(FG)_\tau}b =\Lambda_b(F_\tau,G_\tau).$
Together with the dilation bounds, this gives
\begin{equation}\label{eq:mainregularizedform}
 |\Lambda_{b_\tau}(F,G)|
 \lesssim M\norm F_{A_v^p}\norm G_{A_u^{q\prime}},
 \qquad 0<\tau<1.
\end{equation}

For fixed $\tau<1$, $L_{b_\tau}$ is continuous on $A_\eta^r$ by
Propositions~\ref{prop:finitepair} and~\ref{prop:suppair}:
$b_\tau$ and all its derivatives are bounded, and
$m_j>(1/p-1/q)\gamma_j^+$ when $p\le q$.
For $H\in A_\eta^r$, take its atomic expansion
\[
 H_N=\sum_{k=1}^Nc_k\Phi_k\longrightarrow H,
 \qquad \norm{\mathbf c}_{\ell^r}\lesssim\norm H_{A_\eta^r},
 \qquad \Phi_k=F_kG_k.
\]
If $r\le1$, we can estimate each atom separately. The individual
factor bounds and \eqref{eq:mainregularizedform} give
\begin{align*}
 |L_{b_\tau}(H_N)|
 &\le\sum_{k=1}^N|c_k|\,|\Lambda_{b_\tau}(F_k,G_k)|\\
 &\lesssim M\sum_{k=1}^N|c_k|
 \le M\left(\sum_{k=1}^N|c_k|^r\right)^{1/r}
 \lesssim M\norm H_{A_\eta^r}.
\end{align*}
If $r>1$, we use the finite sign average from the proof of
Theorem~\ref{thm:factor} to retain the coefficient norm in the estimate.
With the factors defined there,
\[
 H_N=2^{-N}\sum_{\boldsymbol\varepsilon\in\{-1,1\}^N}
       F_{\boldsymbol\varepsilon}G_{\boldsymbol\varepsilon},
 \qquad
 \norm{F_{\boldsymbol\varepsilon}}_{A_v^p}
 \norm{G_{\boldsymbol\varepsilon}}_{A_u^{q\prime}}
 \lesssim\left(\sum_{k=1}^N|c_k|^r\right)^{1/r}.
\]
For $q=1$, these are
$F_{\boldsymbol\varepsilon}=\sum_{k=1}^N\varepsilon_kc_kF_k$
and $G_{\boldsymbol\varepsilon}=\sum_{k=1}^N\varepsilon_kG_k$;
the second factor is estimated in $A_u^\infty$.
Applying the form bound to each product in this average gives
\begin{align*}
 |L_{b_\tau}(H_N)|
 &\le2^{-N}\sum_{\boldsymbol\varepsilon}
       |\Lambda_{b_\tau}
          (F_{\boldsymbol\varepsilon},G_{\boldsymbol\varepsilon})|\\
 &\lesssim M\,2^{-N}\sum_{\boldsymbol\varepsilon}
       \norm{F_{\boldsymbol\varepsilon}}_{A_v^p}
       \norm{G_{\boldsymbol\varepsilon}}_{A_u^{q\prime}}\\
 &\lesssim M\left(\sum_{k=1}^N|c_k|^r\right)^{1/r}
 \lesssim M\norm H_{A_\eta^r}.
\end{align*}
In either case, the bound is independent of both $N$ and $\tau$.
As $N\to\infty$, continuity of $L_{b_\tau}$ for fixed
$\tau$ therefore gives, for every $H\in A_\eta^r$,
\begin{equation}\label{eq:mainfunctional}
 |L_{b_\tau}(H)|\lesssim M\norm H_{A_\eta^r},
 \qquad 0<\tau<1.
\end{equation}

When $p>q$, the functional bound gives an integral condition on $b$.
Indeed, $r=s'$ and
$\eta=\bd{\z}{s'\aA}v^{-s'/s}$. Since $s>q$,
\eqref{eq:kernelcondition} implies
$s(a_j+1)>\gamma_j^+-1$ for each $j$.
For each fixed $\tau$, $b_\tau$ is bounded and belongs to $A_v^s$.
The duality in Proposition~\ref{prop:projection}, applied with
$t=s$, therefore gives
\[
 \norm{b_\tau}_{A_v^s}
 \asymp\sup_{\norm H_{A_\eta^{s'}}\le1}
              |\pair H{b_\tau}|
 \lesssim M.
\]
Since $b_\tau(\z)\to b(\z)$, Fatou's lemma now yields
\[
 \int_{\D^n}|b(\z)|^sv(\z)\,dV(\z)
 \le\liminf_{\tau\uparrow1}
       \int_{\D^n}|b(\tau\z)|^sv(\z)\,dV(\z)
 \lesssim M^s.
\]
Consequently $b\in B_s(\om)$ and
$\norm b_{B_s(\om)}\asymp\norm b_{A_v^s}\lesssim M$.

When $p\le q$, the same functional bound gives the derivative
growth condition. To recover $\T b$ at a point, test on the
shifted reproducing kernel. For $\la\in\D^n$, take
$H_{\la}(\w)=\prod_j (1-\overline{\lambda_j}w_j)^{-a_j-m_j-2}.$
The kernel estimate shows that $H_{\la}\in A_\eta^r$ and
\begin{align*}
 \norm{H_{\la}}_{A_\eta^r}
 &\lesssim[\eta(\la)\bd{\la}{2}]^{1/r}
           \prod_j(1-|\lambda_j|^2)^{-a_j-m_j-2}\\
 &=\ow{\la}^d\bd{\la}{-\m}.
\end{align*}
Here $m_j>d\gamma_j^+$ ensures convergence of the kernel integral,
and the last equality is \eqref{eq:massidentity}.
Using $H_{\la}$ in \eqref{eq:mainfunctional} and the reproducing
identity \eqref{eq:kernelpair}, we obtain
\[
|\T b_\tau(\la)|=|\pair{H_{\la}}{b_\tau}| \lesssim M\ow{\la}^d\bd{\la}{-\m}, \qquad 0<\tau<1.
\]
At each fixed $\la$, holomorphic convergence gives
$\T b_\tau(\la)\to\T b(\la)$. Let $\tau\uparrow1$,
multiply by $\bd{\la}{\m}\ow{\la}^{-d}$ and take the supremum
over $\la\in\D^n$ to conclude that
$\mathcal N_{p,q}(b)\lesssim M$.
Together with the sufficiency estimates and the equivalence of
Besov and Bergman norms, these bounds prove both assertions.
\end{proof}

\begin{proof}[Proof of Theorem~\ref{thm:main} for $0<q<1$]
As in the preceding case, \eqref{eq:model} allows us to work with
$\Hb:A_v^p\to\overline{A_v^q}$.

\emph{Sufficiency.}
Let $F\in A_v^p$. We apply \eqref{eq:outputestimate}, verifying
the finiteness of its right-hand side in each range.
Suppose $q<p$ and set $s=pq/(p-q)$.
Since $q/p+q/s=1$, H\"older's inequality gives
\begin{align*}
 \norm{\Hb F}_{A_v^q}^q
 &\lesssim\int|F|^q|\T b|^q\bd{\z}{q\m}v\,dV\\
 &\le\left(\int|F|^pv\,dV\right)^{q/p}
       \left(\int|\T b|^s\bd{\z}{s\m}v\,dV\right)^{q/s}\\
 &\lesssim\norm F_{A_v^p}^q\norm b_{A_v^s}^q.
\end{align*}
In the last step we used \eqref{eq:LP} with exponent $s$;
this identity remains valid when $s<1$.

Suppose instead that $p\le q$ and $\mathcal N_{p,q}(b)<\infty$.
Here the point evaluation estimate \eqref{eq:point} gives
$|F(\z)|\ow{\z}^{1/p}\lesssim\norm F_{A_v^p}.$
After inserting the symbol bound in \eqref{eq:outputestimate},
we use this pointwise estimate to control the extra power of the input:
\begin{align*}
 \norm{\Hb F}_{A_v^q}^q
 &\lesssim\mathcal N_{p,q}(b)^q
       \int|F|^q\ow{\z}^{q/p-1}v\,dV\\
 &=\mathcal N_{p,q}(b)^q\int
       |F|^pv\big(|F|\ow{\z}^{1/p}\big)^{q-p}\,dV\\
 &\lesssim\mathcal N_{p,q}(b)^q
       \norm F_{A_v^p}^{q-p}\norm{F}_{A_v^p}^p
 \le C\mathcal N_{p,q}(b)^q\norm F_{A_v^p}^q.
\end{align*}
For $p=q$ the middle factor has exponent zero and is simply omitted.

In both ranges the right-hand side of \eqref{eq:outputestimate}
is finite. The shifted integral therefore defines the bounded
extension to every $F\in A_v^p$, with
\[
 \norm{\Hb}_{A_v^p\to\overline{A_v^q}}
 \lesssim
 \begin{cases}
 \norm b_{A_v^s},&q<p,\\
 \mathcal N_{p,q}(b),&p\le q.
 \end{cases}
\]

\emph{Necessity.}
Write $M=\norm{\Hb}_{A_v^p\to\overline{A_v^q}}$.
For any $\la\in\D^n$, choose the normalized test function
\[
 F_{\la}(\w)=\ow{\la}^{-1/p}
       \prod_j\frac{(1-|\lambda_j|^2)^{m_j}}
                      {(1-\overline{\lambda_j}w_j)^{m_j}}.
\]
The kernel estimate gives $\norm{F_{\la}}_{A_v^p}\lesssim1$.
The exponent in this test function is chosen so that evaluation
of its image at $\la$ combines the two kernel factors into the
kernel for the shifted derivative:
\begin{align*}
 (\Hb F_{\la})(\la)
 &=\ow{\la}^{-1/p}\bd{\la}{\m}
     \int\prod_j(1-\overline{\lambda_j}w_j)^{-a_j-m_j-2}
                    \overline{b(\w)}\,d\nu_{\aA}(\w)\\
 &=\ow{\la}^{-1/p}\bd{\la}{\m}
                       \overline{\T b(\la)}.
\end{align*}
The final equality is \eqref{eq:kernelpair}. This identity uses
only the initial integrability of $b$.

If $p\le q$, boundedness of point evaluation in the output space gives
\begin{align*}
 \bd{\la}{\m}|\T b(\la)|\ow{\la}^{1/q-1/p}
 &=\ow{\la}^{1/q}|(\Hb F_{\la})(\la)|\\
 &\lesssim\norm{\Hb F_{\la}}_{A_v^q}
 \le M\norm{F_{\la}}_{A_v^p}\lesssim M.
\end{align*}
Taking the supremum proves \eqref{eq:symbolsup} with the required
norm bound; a single test at each point suffices in this range.

Now suppose $q<p$. To obtain the integral condition on the symbol,
we combine the tests over all cells. Using the points $\la_k\in\overline{Q_k}$ and
the functions $F_k$ from \eqref{eq:atoms}, we substitute the preceding
evaluation identity into \eqref{eq:randomtest} to obtain
\begin{equation}\label{eq:diagonal}
 \sum_k|c_k|^q
   \left[\bd{\la_k}{\m}|\T b(\la_k)|
                  \ow{\la_k}^{1/q-1/p}\right]^q
 \lesssim M^q\left(\sum_k|c_k|^p\right)^{q/p}.
\end{equation}
Write $A_k$ for the nonnegative bracket in \eqref{eq:diagonal}.
We choose the coefficients to recover its $\ell^s$ norm.
On a finite set of indices let $S=\sum_kA_k^s$, where
$s=pq/(p-q)$. If $S>0$, choose
$c_k=A_k^{s/p}S^{-1/p}$. Then
\[
 \sum_k|c_k|^p=1,\qquad
 \sum_k|c_k|^qA_k^q
 =S^{-q/p}\sum_kA_k^{sq/p+q}
 =S^{1-q/p}=S^{q/s}.
\]
Thus $S^{1/s}\lesssim M$; the assertion is immediate when $S=0$.
Since $s(1/q-1/p)=1$, this proves for every finite set
\begin{equation}\label{eq:discretesymbol}
 \sum_k|\T b(\la_k)|^s\bd{\la_k}{s\m}\ow{\la_k}
 \lesssim M^s.
\end{equation}
Since the discrete bound is uniform in the representatives, choose
$\la_k$ at a maximum of $|\T b|$ on
$\overline{Q_k}$. Local comparability gives
\begin{align*}
 \int_{Q_k}|\T b(\z)|^s\bd{\z}{s\m}v(\z)\,dV(\z)
 &\lesssim|\T b(\la_k)|^s\bd{\la_k}{s\m}
                           \int_{Q_k}v\,dV\\
 &\lesssim|\T b(\la_k)|^s\bd{\la_k}{s\m}\ow{\la_k}.
\end{align*}
Summing first over finite collections of cells and then exhausting
the polydisc, we obtain from \eqref{eq:discretesymbol} and monotone
convergence that
\[
 \int_{\D^n}|\T b|^s\bd{\z}{s\m}v\,dV\lesssim M^s.
\]
The derivative norm identity \eqref{eq:LP} now gives the required bound
$\norm b_{A_v^s}\lesssim M$.
Returning to the Besov norms through \eqref{eq:model} completes
the necessity estimates and the asserted norm equivalences.
\end{proof}

Moment derivatives also describe the one-variable Bloch space; see
\cite[Theorem~3]{MorenoPelaezRosa}. For the polydisc, the
coordinate-face norm in \eqref{eq:Blochdef} is equivalent to
Constantin's norm on $\mathcal B^{(n)}$; see
\cite[Section~1, pp.~160--161]{ConstantinWeak}.
The higher-order product-derivative characterization of this space is
given in \cite[Section~3.2, (9), p.~5]{Sehba}.
We only need to adapt its derivative operator to the parameters
$\aA$ of our kernel.

\begin{lemma}\label{lem:bloch}
For $a_j>-1$ and positive integers $m_j$,
\begin{equation}\label{eq:bloch}
 \norm b_{\Bp}\asymp\sup_{\z\in\D^n}
          \bd{\z}{\m}|\T b(\z)|.
\end{equation}
\end{lemma}
\begin{proof}
Sehba's derivative of order $\m$ is
$\prod_j\prod_{\ell=0}^{m_j-1}(R_j+2+\ell)$, a fixed scalar
multiple of $T_{\mathbf0,\m}$. Thus the cited characterization gives
\[
 \norm b_{\Bp}\asymp
 \sup_{\z\in\D^n}\bd{\z}{\m}
                    |T_{\mathbf0,\m}b(\z)|.
\]
The separate term $|b(\mathbf0)|$ in that characterization is
already controlled by this supremum.
On the holomorphic growth space with norm
$\norm g_{\m,\infty}=\sup\bd{\z}{\m}|g(\z)|$, the resolvent
\[
 (R_j+c)^{-1}g(\z)
 =\int_0^1t^{c-1}
       g(z_1,\ldots,tz_j,\ldots,z_n)\,dt,
 \qquad c>0,
\]
has norm at most $1/c$, since $1-|z_j|^2\le1-t^2|z_j|^2$.
Consequently, both
\[
 \frac{R_j+a_j+2+\ell}{R_j+2+\ell}
 =I+a_j(R_j+2+\ell)^{-1}
 \quad\text{and}\quad
 \frac{R_j+2+\ell}{R_j+a_j+2+\ell}
 =I-a_j(R_j+a_j+2+\ell)^{-1}
\]
are bounded on this growth space. Multiplying these finitely many
commuting operators shows that
$\norm{T_{\aA,\m}b}_{\m,\infty}
\asymp\norm{T_{\mathbf0,\m}b}_{\m,\infty}$, proving the assertion.
\end{proof}

\begin{corollary}\label{cor:diagonal}
For $p=q>0$ satisfying \eqref{eq:kernelcondition},
$\Hb:B_p(\om)\longrightarrow\overline{B_p(\om)} \text{ is bounded}\quad\Longleftrightarrow\quad b\in\Bp,$
with comparable norms. Every $b\in\Bp$ belongs to
$L^1(\nu_{\aA})$.
\end{corollary}
\begin{proof}
For $p=q$, the weight factor in $\mathcal N_{p,p}$ is one.
Theorem~\ref{thm:main} and Lemma~\ref{lem:bloch} therefore give the
boundedness criterion and norm equivalence. By
\cite[Lemma~3.3, p.~6]{Sehba},
\[
 |b(\z)|\lesssim\norm b_{\Bp}
       \prod_j\log\frac4{1-|z_j|^2}.
\]
Each logarithmic factor is integrable against
$(1-|z_j|^2)^{a_j}dA(z_j)$ for $a_j>-1$, so
$b\in L^1(\nu_{\aA})$.
\end{proof}

\begin{remark}
The kernel in \cite[Theorems~6--7 of the preprint]{Wang} is induced by the target
weight; here it is prescribed independently by $\nu_{\aA}$.
This distinction matters when comparing the ranges of the two results.
For example, $v=1$ and $a=0$ in one variable give
$\gamma^+=2$, so the condition $q(a+2)>\gamma^+$ fails for $q<1$.
The present theorem therefore does not include that matched-kernel case.
The interval $\gamma_j^++q-1<q(a_j+2)\le\gamma_j^+$ also lies
outside the range established by our sufficiency proof for $q<1$.
\end{remark}

\end{document}